\documentclass[aap,MSNbibl,citesort,dvips]{arximspdf}
\usepackage{algpseudocode}
\usepackage{algorithm,algorithmicx}
\usepackage{mathbh}

\doi{10.1214/10-AAP735}
\volume{21}
\issue{6}
\pubyear{2011}
\firstpage{2109}
\lastpage{2145}

\makeatletter

\let\epsilon\varepsilon

\newtheorem{prop}[thm]{Proposition}
\newtheorem{lem}[thm]{Lemma}
\newtheorem{cor}[thm]{Corollary}

\newproclaim{assum}{Assumption}
\newproclaim{rem}{Remark}

\newcommand{\iint}{\int\!\!\!\int}
\newcommand{\idotsint}{\int\cdots\int}
\newcommand{\mcf}[2]{\mathcal{F}_{#1}^{#2}}

\newcommand{\eqref}[1]{(\ref{#1})}

\newcommand{\wrt}{with respect to}
\newcommand{\iid}{i.i.d.}
\newcommand{\as}{a.s.}

\newcommand{\rset}{\mathbb{R}}
\newcommand{\nset}{\mathbb{N}}
\newcommand{\lone}{\mathsf{L}^1}

\newcommand{\1}{\mathbf{1}}
\newcommand{\eqdef}{\stackrel{\mathrm{def}}{=}}
\newcommand{\rmd}{d}
\newcommand{\rme}{e}
\newcommand{\mcb}[1]{\mathcal{F}_{\mathrm{b}}(#1)}
\newcommand{\supnorm}[1]{|#1|_{\infty}}
\newcommand{\esssup}[2][]%
{\ifthenelse{\equal{#1}{}}{| #2 |_\infty}{| #2|^2_{\infty}}}
\newcommand{\oscnorm}[2][]%
{\ifthenelse{\equal{#1}{}}{\ensuremath{\operatorname{osc}(#2
)}}{\ensuremath{\operatorname{osc}^{#1}(#2)}}}

\newcommand{\PP}{\mathbb{P}}
\newcommand{\PE}{\mathbb{E}}
\newcommand{\CPE}[3][]
{\ifthenelse{\equal{#1}{}}{\mathbb{E}[ #2  | #3
]}{\mathbb{E}_{#1}[ #2  | #3 ]}}
\newcommand{\CPP}[3][]
{\ifthenelse{\equal{#1}{}}{\mathbb{P}[ #2  | #3
]}{\mathbb{P}_{#1}[ #2  | #3 ]}}

\newcommand{\dlim}{\stackrel{\mathcal{D}}{\longrightarrow}}
\newcommand{\plim}{\stackrel{\mathrm{P}}{\longrightarrow}}

\newcommand{\Xset}{\mathbb{X}}
\newcommand{\Xsigma}[1][]%
{%
\ifthenelse{\equal{#1}{}}{\ensuremath{\mathcal{B}(\Xset)}}{\ensuremath{
\mathcal{B}(\Xset^{#1})}}
}
\newcommand{\Yset}{\mathbb{Y}}
\newcommand{\Ysigma}{\mathcal{B}(\Yset)}
\newcommand{\chunk}[4][]%
{\ifthenelse{\equal{#1}{}}{\ensuremath{{#2}_{#3:#4}}}{
\ensuremath{#2^#1}_{#3:#4}}
}

\newcommand{\M}{M}
\newcommand{\Xinit}{\chi}
\newcommand{\XinitIS}[2][]
{\ifthenelse{\equal{#1}{}}{\ensuremath{\rho_{#2}}}{\ensuremath{\check{
\rho}_{#2}}}}

\newcommand{\adjfunc}[4][]
{\ifthenelse{\equal{#1}{}}{\ifthenelse{\equal{#4}{}}{\vartheta_{#2}}{
\vartheta_{#2}(#4)}}
{\ifthenelse{\equal{#1}{smooth}}{\ifthenelse{\equal{#4}{}}{\tilde{
\vartheta}_{#2}}{\tilde{\vartheta}_{#2}(#4)}}
{\ifthenelse{\equal{#1}{fully}}{\ifthenelse{\equal{#4}{}}{\vartheta^
\star_{#2}}{\vartheta^\star_{#2}(#4)}}{\mathrm{erreur}}}}}

\newcommand{\post}[3][]%
{
\ifthenelse{\equal{#1}{}}{\ensuremath{\phi_{#2|#3}}}%
{\ifthenelse{\equal{#1}{hat}}{\ensuremath{\phi^{N}_{#2|#3}}}
{\ifthenelse{\equal{#1}{tilde}}{\ensuremath{\tilde{\phi}^{N}_{#2|#3}}}
{\ifthenelse{\equal{#1}{tar}}{\ensuremath{\phi^{N,
\mathrm{t}}_{#2|#3}}}}
}
}
}

\newcommand{\filt}[2][]%
{
\ifthenelse{\equal{#1}{}}{\ensuremath{\phi_{#2}}}%
{\ifthenelse{\equal{#1}{hat}}{\ensuremath{\phi^{N}_{#2}}}
{\ifthenelse{\equal{#1}{tilde}}{\ensuremath{\tilde{\phi}^{N}_{#2}}}
{\ifthenelse{\equal{#1}{tar}}{\ensuremath{\phi^{N,\mathrm{t}}_{#2}}}
{\ifthenelse{\equal{#1}{aux}}{\ensuremath{\phi^{N,\mathrm{a}}_{#2}}}
}
}
}
}
}

\newcommand{\unfilt}[2][]%
{
\ifthenelse{\equal{#1}{}}{\ensuremath{\gamma_{#2}}}%
{\ifthenelse{\equal{#1}{hat}}{\ensuremath{\gamma^{N}_{#2}}}
{\ifthenelse{\equal{#1}{tilde}}{\ensuremath{\tilde{\gamma}^{N}_{#2}}}
{\ifthenelse{\equal{#1}{tar}}{\ensuremath{\gamma^{N,\mathrm{t}}_{#2}}}
{\ifthenelse{\equal{#1}{aux}}{\ensuremath{\gamma^{N,\mathrm{a}}_{#2}}}
}
}
}
}
}

\newcommand{\asymVar}[4][]{
\ifthenelse{\equal{#1}{}}{\ifthenelse{\equal{#4}{}}{\ensuremath{
\Gamma_{#2|#3}}}{\ensuremath{\Gamma_{#2|#3}[#4]}}}
{\ifthenelse{\equal{#4}{}}{\ensuremath{\Gamma_{#1,#2|#3}}}{\ensuremath{
\Gamma_{#1,#2|#3}[#4]}}}
}

\newcommand{\incrasymVar}[4][]{
\ifthenelse{\equal{#1}{}}{\ensuremath{\sigma^2_{#2,#3}[#4
]}}{\ensuremath{\sigma^2_{#1,#2,#3}[#4]}}
}

\newcommand{\BK}[1]{\mathrm{B}_{#1}}

\newcommand{\kiss}[3][]
{\ifthenelse{\equal{#1}{}}{p_{#2}}
{\ifthenelse{\equal{#1}{fully}}{p^{\star}_{#2}}
{\ifthenelse{\equal{#1}{smooth}}{\tilde{r}_{#2}}{\mathrm{erreur}}}}}

\newcommand{\Kiss}[3][]
{\ifthenelse{\equal{#1}{}}{P_{#2}}
{\ifthenelse{\equal{#1}{fully}}{P^{\star}_{#2}}
{\ifthenelse{\equal{#1}{smooth}}{\tilde{R}_{#2}}{\mathrm{erreur}}}}}

\newcommand{\epart}[2]{\ensuremath{\xi_{#1}^{#2}}}
\newcommand{\ewght}[2]{\ensuremath{\omega_{#1}^{#2}}}
\newcommand{\ewghtfunc}[1]{\ensuremath{\omega_{#1}}}
\newcommand{\sumwght}[2][]{%
\ifthenelse{\equal{#1}{}}{\ensuremath{\Omega_{#2}}}{\ensuremath{
\Omega_{#2}^{(#1)}}}}
\newcommand{\instrpostaux}[2]{\pi_{#1|#2}}

\newcommand{\F}{\mathcal L}
\newcommand{\extens}[3]{{\Pi_{#2,#3} #1}}

\makeatother

\begin{document}
\begin{frontmatter}

\title{Sequential {M}onte {C}arlo smoothing for general state space
hidden Markov models\thanksref{TITL1}}
\runtitle{SMC smoothing for HMM}
\thankstext{TITL1}{Supported in part by the French National Research
Agency, under the program ANR-07 ROBO 0002 CFLAM.}

\begin{aug}
\author[A]{\fnms{Randal} \snm{Douc}\ead[label=e1]{randal.douc@it-sudparis.eu}},
\author[B]{\fnms{Aur\'elien} \snm{Garivier}\ead[label=e2]{aurelien.garivier@telecom-paristech.fr}},
\author[B]{\fnms{Eric}
\snm{Moulines}\corref{}\ead[label=e3]{eric.moulines@telecom-paristech.fr}}\break
and
\author[C]{\fnms{Jimmy}~\snm{Olsson}\ead[label=e4]{jimmy@maths.lth.se}}
\runauthor{Douc, Garivier, Moulines and Olsson}
\affiliation{Institut T\'el\'ecom/T\'el\'ecom SudParis, CNRS UMR 5157,
CNRS/T\'el\'ecom ParisTech, CNRS UMR 5141, Institut T\'el\'ecom/T\'el\'
ecom ParisTech, CNRS UMR 5141 and~Lund University}
\address[A]{R. Douc\\
T\'el\'ecom SudParis \\
9 rue Charles Fourier \\
91011 Evry-Cedex\\
France \\
\printead{e1}}
\address[B]{A. Garivier\\
E. Moulines\\
T\'el\'ecom ParisTech\\
46 rue Barrault \\
75634 Paris, Cedex 13\\
France \\
\printead{e2}\\
\hphantom{E-mail: }\printead*{e3}}
\address[C]{J. Olsson\\
Center for Mathematical Sciences \\
Lund University \\
Box 118\\
SE-22100 Lund\\
Sweden \\
\printead{e4}}
\end{aug}

\received{\smonth{7} \syear{2009}}
\revised{\smonth{3} \syear{2010}}

%
\begin{abstract}
Computing smoothing distributions, the distributions of one
or more states conditional on past, present, and future observations is
a recurring problem when operating on general hidden Markov models. The
aim of this paper is to provide a foundation of particle-based
approximation of such distributions and to analyze, in a common
unifying framework, different schemes producing such approximations.
In this setting, general convergence results, including exponential
deviation inequalities and central limit theorems, are established. In
particular, time uniform bounds on the marginal smoothing error are
obtained under appropriate mixing conditions on the transition kernel
of the latent chain. In addition, we propose an algorithm approximating
the joint smoothing distribution at a cost that grows only linearly
with the number of particles.
\end{abstract}

%

\end{frontmatter}

\section{Introduction}\label{sec1}

Statistical inference in general state space\vspace*{-2pt} hidden Mar\-kov models (HMM)
involves computation of the \emph{posterior distribution} of a~set
$\chunk{X}{s}{s'} \eqdef[X_s, \dots, X_{s'}]$ of state variables
conditional on a record $\chunk{Y}{0}{T} = \chunk{y}{0}{t}$ of
observations. This distribution will, in the following, be denoted by
$\post{s:s'}{T}$ where the dependence of this measure on the observed
values $\chunk{y}{0}{T}$ is implicit. The posterior distribution can
be expressed in closed-form only in very specific cases, principally,
when the state space model is linear and Gaussian or when the state
space of the hidden Markov chain is a finite set. In the vast majority
of cases, nonlinearity or non-Gaussianity render analytic solutions
intractable \cite
{kailathsayedhassibi2000,risticarulampalamgordon2004,cappemoulinesryden2005,vanhandel2007}.\vadjust{\goodbreak}

This limitation has led to an increase of interest in alternative
computational strategies handling more general state and measurement
equations without constraining a priori the behavior of the posterior
distributions. Among these, \emph{sequential {M}onte {C}arlo} (SMC)
\emph{methods} play a central role. SMC methods---in which the \emph
{sequential importance sampling} and \emph{sampling importance
resampling} methods proposed by \cite{handschinmayne1969} and \cite
{rubin1987}, respectively, are combined---refer to a class of algorithms
approximating a \emph{sequence of probability distributions}, defined
on a \emph{sequence of probability spaces}, by updating recursively a
set of random \emph{particles} with associated nonnegative \emph
{importance weights}. The SMC methodology has emerged as a key tool for
approximating state posterior distribution flows in general state space
models; see \cite{delmoralmiclo2000,delmoral2004,delmoraldoucet2009}
for general introductions as well as theoretical results for SMC
methods and \cite{liu2001,doucetdefreitasgordon2001,risticarulampalamgordon2004} for
applications of SMC within a variety of scientific fields.

The recursive formulas generating the \emph{filter distributions}
$\phi_{T}$ (short-hand notation for $\post{T:T}{T}$) and the \emph
{joint smoothing distributions} $\post{0:T}{T}$ are closely related;
thus, executing the standard SMC scheme in the filtering mode provides,
as a by-product, approximations of the joint smoothing distributions.
More specifically, the branches of the genealogical tree associated
with the historical evolution of the filtering particles up to time
step $T$ form, when combined with the corresponding importance weights
of these filtering particles, a weighted sample approximating the joint
smoothing distribution $\post{0:T}{T}$; see \cite{delmoral2004}, Section 3.4,  for details. From these paths, one may readily
obtain a weighted sample targeting the fixed lag or fixed interval
smoothing distribution by extracting the required subsequence of states
while retaining the weights. This appealingly simple scheme can be used
successfully for estimating the joint smoothing distribution for small
values of $T$ or any marginal smoothing distribution $\post{s}{T}$,
with $s \leq T$, when $s$ and $T$ are close; however, when $T$ is large
and $s \ll T$, the associated particle approximations are inaccurate
since the genealogical tree degenerates gradually as the interacting
particle system evolves \cite{godsilldoucetwest2004,fearnheadwyncolltawn2010}.

In this article, we thus give attention to more sophisticated
approaches and consider instead the \emph{forward filtering backward
smoothing} (FFBSm) \emph{algorithm} and the \emph{forward filtering
backward simulation} (FFBSi) \emph{sampler}. These algorithms share
some similarities with the Baum--Welch algorithm for finite state space
models and the Kalman filter-based smoother and simulation smoother for
linear Gaussian state space models \cite{dejongshephard1995}. In the
FFBSm algorithm, the particle weights obtained when approximating the
filter distributions in a forward filtering pass are modified in a
backward pass; see \cite{kitagawa1996,huerzelerkuensch1998,doucetgodsillandrieu2000}. The FFBSi algorithm simulates, conditionally
independently given the particles and particle weights produced in a
similar forward filtering pass, state trajectories being approximately
distributed according to the joint smoothing distribution; see \cite
{godsilldoucetwest2004}.

The computational complexity of the FFBSm algorithm when used for
estimating marginal fixed interval smoothing distributions or of the
original formulation of the FFBSi sampler grows (in most situations) as
the square of the number $N$ of particles multiplied by the time
horizon $T$. To alleviate this potentially very large computational
cost, some methods using intricate data structures for storing the
particles have been developed; see, for example,~\cite
{klaasbriersdefreitasdoucetmaskelllang2006}. These algorithms have a
complexity of order $O(N \log(N))$ and are thus amenable to practical
applications; however, this reduction in complexity comes at the cost
of introducing some level of approximation.

In this paper, a modification of the original FFBSi algorithm is
presented. The proposed scheme has a complexity that grows only \emph
{linearly} in $N$ and does not involve any numerical approximation
techniques. This algorithm may be seen as an alternative to a recent
proposal by \cite{fearnheadwyncolltawn2010} which is based on the
so-called \emph{two-filter algorithm} \cite{briersdoucetmaskell2010}.

\setcounter{footnote}{1}
The smoothing weights computed in the backward pass of the FFBSm
algorithm at a given time instant $s$ (or the law of the FFBSi
algorithm) are statistically dependent on all forward filtering pass
particles and weights computed before and after this time instant. This
intricate dependence structure makes the analysis of the resulting
particle approximation challenging; up to our best knowledge, only a
single consistency result is available in~\cite
{godsilldoucetwest2004}, but its proof is plagued by a (subtle) mistake
that seems difficult to correct. Therefore, very little is known about
the convergence of the schemes under consideration, and the second
purpose of this paper is to fill this gap.\footnote{Since the first version of this paper has been released, an
article \cite{delmoraldoucetsingh2010} has been published. This work,
developed completely independently from ours, complement the results
presented in this manuscript. In particular, this paper presents a
functional central limit theorems as well as nonasymptotic variance
bounds. Additionally, this work shows how the forward filtering
backward smoothing estimates of additive functionals can be computed
using a forward only recursion.}
In this contribution, we focus first on finite time horizon
approximations. Given a finite time horizon $T$, we derive \emph
{exponential deviation inequalities} stating that the probability of
obtaining, when replacing $\post{s:T}{T}$ by the corresponding FFBSm
or FFBSi estimator, a {M}onte {C}arlo error exceeding a~given $\epsilon
> 0$ is bounded by a quantity of order $O(\exp(- c N \epsilon^2))$
where $c$ is positive constant depending on $T$ as well as the target
function under consideration. The obtained inequalities, which are
presented in Theorem \ref{thm:Hoeffding-FFBS} (FFBSm) and Corollary~\ref{cor:Hoeffding-FFBSi} (FFBSi), hold for any given number $N$ of
particles and are obtained by combining a novel backward error
decomposition with an adaptation of the Hoeffding inequality to
statistics expressed as ratios of random variables. We then consider
the asymptotic (as the number~$N$ of particles tends to infinity)
regime and establish a central limit theorem (CLT) with rate $\sqrt
{N}$ and with an explicit expression of the asymptotic variance; see
Theorem \ref{thm:FFBS-CLT}. The proof of our CLT relies on a
technique, developed gradually in \cite
{chopin2004,kuensch2005,doucmoulines2008}, which is based on a CLT for
triangular arrays of dependent random variables; however, since we are
required to take the complex dependence structure of the smoothing
weights into account, our proof is significantly more involved than in
the standard filtering framework considered in the mentioned works.

The second part of the paper is devoted to time uniform results, and we
here study the behavior of the particle-based marginal smoothing
distribution approximations as the time horizon $T$ tends to infinity.
In this setting, we first establish, under the assumption that the
Markov transition kernel~$\M$ of the latent signal is strongly mixing
(Assumption \ref{assum:strong-mixing-condition}), time uniform deviation bounds of the type described
above which hold for any particle population size $N$ and where the
constant $c$ is \emph{independent} of $T$; see Theorem \ref
{theo:Hoeffding-uniform}. This result may seem surprising, and the
nonobvious reason for its validity stems from the fact that the
underlying Markov chain forgets, when evolving conditionally on the
observations, its initial conditions in the forward \emph{as well as}
the backward directions. Finally, we prove (see Theorem \ref
{theo:CLT-uniform}), under the same uniform mixing assumption, that the
asymptotic variance of the CLT for the particle-based marginal
smoothing distribution approximations remains bounded as $T$ tends to
infinity. The uniform mixing assumption in Assumption \ref{assum:strong-mixing-condition} points
typically to applications where the state space of the latent signal is
compact; nevertheless, in the light of recent results on filtering
stability \cite
{kleptsynaveretennikov2008,doucfortmoulinespriouret2009} one may
expect the geometrical contraction of the backward kernel to hold for a
significantly larger class of nonuniformly mixing models (see \cite
{doucfortmoulinespriouret2009} for examples from, e.g.,
financial economics). But even though the geometrical mixing rate is
supposed to be constant in this more general case, applying the
mentioned results will yield a bound of contraction containing a
multiplicative constant depending highly on the initial distributions
as well as the observation record under consideration. Since there are
currently no available results describing this dependence, applying
such bounds to the instrumental decomposition used in the proof of
Theorem \ref{thm:Hoeffding-FFBS} seems technically involved. Recently,
\cite{vanhandel2009} managed to derive \emph{qualitative} time
average convergence results for standard (bootstrap-type) particle
filters under a mild tightness assumption being satisfied also in the
noncompact case when the hidden chain is geometrically ergodic. Even
though this technique does not (on the contrary to our approach) supply
a rate of convergence, it could possibly be adopted to our framework in
order to establish time average convergence of the particle-based
marginal smoothing distribution approximations in a noncompact setting.

The paper is organized as follows. In Section \ref{sec:FFBS}, the
FFBSm algorithm and the FFBSi sampler are introduced. An exponential
deviation inequality for the fixed interval joint smoothing
distribution is derived in Section \ref{sec:exponentialFFBS}, and a
CLT is established in Section \ref{sec:CLTFFBS}. In Section \ref
{sec:TimeUniformExponentialFFBS}, time uniform exponential bounds on
the error of the FFBSm marginal smoothing distribution estimator are
computed under the mentioned mixing condition on the kernel $\M$.
Finally, under the same mixing condition, an explicit bound on the
asymptotic variance of the marginal smoothing distribution estimator is
derived in Section \ref{sec:TimeUniformCLTFFBS}.

\subsection*{Notation and definitions}
\label{sec:Notations and Definitions}
For any sequence $\{a_n\}_{n \geq0}$ and\vspace*{-2pt} any pair of integers $0 \leq
m \leq n$, we denote $\chunk{a}{m}{n} \eqdef(a_m, \dots, a_n)$.
We assume in the following that all random variables are defined on a
common probability space $(\Omega, \mcf{}{}, \PP)$. The sets $\Xset
$ and $\Yset$
are supposed to be Polish spaces and we denote by $\Xsigma$ and
$\Ysigma$ the associated Borel $\sigma$-algebras. $\mcb{\Xset}$
denotes the set of all bounded $\Xsigma/ \mathcal{B}(\rset
)$-measurable functions from $\Xset$ to $\rset$. For any\vspace*{-2pt}
measure~$\zeta$ on $(\Xset, \Xsigma)$ and any $\zeta$-integrable function
$f$, we set $\zeta(f) \eqdef\int_\Xset f(x) \zeta(\rmd x)$. Two
measures $\zeta$ and $\zeta'$ are said to be \emph{proportional}
(written $\zeta\propto\zeta'$) if they differ only by a
normalization constant.

A kernel $V$ from $(\Xset, \Xsigma)$ to $(\Yset, \Ysigma)$ is a
mapping from $\Xset\times\Ysigma$ into $[0, 1]$ such that, for each
$A \in\Ysigma$, $x \mapsto V (x, A)$ is a nonnegative, bounded, and
measurable function on $\Xset$, and, for each $x \in\Xset$, $A
\mapsto V (x, A)$\vspace*{-1pt} is a measure on $\Ysigma$. For $f \in\mcb{\Xset}$
and $x \in\Xset$, denote by $V(x, f) \eqdef\int V(x, \rmd x')
f(x')$; we will sometimes also use the abridged notation $Vf(x)$
instead of $V(x, f)$. For a measure $\nu$ on $(\Xset, \Xsigma)$, we
denote by $\nu V$ the measure on $(\Yset, \Ysigma)$ defined by, for\vspace*{-1pt}
any $A \in\Ysigma$, $\nu V(A) \eqdef\int_\Xset V(x,A) \nu(\rmd x)$.

Consider now a possibly nonlinear state space model, where the \emph
{state process} $\{ X_t \}_{t \geq0}$ is a Markov chain on the state
space $(\Xset, \Xsigma)$. Even though~$t$ is not necessarily a
temporal index, we will often refer to this index as ``time.'' We
denote by $\Xinit$ and $\M$ the initial distribution and transition
kernel, respectively, of this process. The state process is assumed to be
hidden but partially observed through the \emph{observations} $\{ Y_t
\}_{t \geq0}$ which are $\Yset$-valued random variables being
conditionally independent given the latent state sequence $\{ X_t \}_{t
\geq0}$; in addition, there exists a $\sigma$-finite measure $\lambda
$ on $(\Yset, \Ysigma)$ and a nonnegative transition density function
$g$ on $\Xset\times\Yset$ such that $\CPP{Y_t \in A}{X_t} = \int_A
g(X_t, y) \lambda(\rmd y)$ for all $A \in\Ysigma$. The mapping $x
\mapsto g(x, y)$ is referred to as the \emph{likelihood function} of
the state given an observed value $y \in\Yset$. The kernel $\M$ as
well as the transition density $g$ are supposed to be known. In the
setting of this paper, we assume that we have access to a record of\vspace*{-2pt}
arbitrary but fixed observations $\chunk{y}{0}{T} \eqdef[y_0, \dots,
y_T]$, and our main task is to estimate the posterior distribution of
(different subsets of)\vspace*{-1pt} the state vector $X_{0:T}$ given these
observations. For any $t \geq0$, we denote by $g_t(x) \eqdef g(x,
y_t)$ (where the dependence on $y_t$ is implicit) the likelihood
function of the state $X_t$ given the observation $y_t$.

For simplicity, we consider a \emph{fully dominated} state space model
for which there exists a $\sigma$-finite measure $\nu$ on $(\Xset,
\Xsigma)$ such that, for all $x \in\Xset$, $\M(x, \cdot)$ has a
transition probability density $\ensuremath{m}(x, \cdot)$ \wrt\ $\nu
$. For
notational simplicity, $\nu(\rmd x)$ will sometimes be replaced by
$\rmd x$.

For any initial distribution $\Xinit$ on $(\Xset, \Xsigma)$ and any
$0 \leq s \leq s'\leq T$, denote by $\post{s:s'}{T}$ the posterior
distribution of the state vector $\chunk{X}{s}{s'}$ given the
observations $\chunk{y}{0}{T}$. For lucidity, the dependence of $\post
{s:s'}{T}$ on the initial distribution~$\Xinit$ is omitted. Assuming
that $\idotsint\Xinit(\rmd x_0) \prod_{u = 1}^T g_{u-1}(x_{u-1}) \M
(x_{u-1},\allowbreak \rmd x_u) g_T(x_T)\,{>}\,0$, this distribution may be expressed
as, for all $h\,{\in}\,\mcb{\Xset^{s'-s+1}}$,
\[
\post{s:s'}{T}(h) =
\frac{\idotsint\Xinit(\rmd x_0) \prod_{u = 1}^T g_{u-1}(x_{u-1}) \M
(x_{u-1}, \rmd x_u) g_T(x_T) h(\chunk{x}{s}{s'})}
{\idotsint\Xinit(\rmd x_0) \prod_{v = 1}^T g_{v-1}(x_{v-1}) \M
(x_{v-1}, \rmd x_v) g_T(x_T)} .
\]
In the expression above, the dependence on the observation sequence is
implicit. If $s = s'$, we use $\post{s}{T}$ (the marginal smoothing
distribution at\vspace*{-2pt} time~$s$) as shorthand for $\post{s:s}{T}$. If $s = s'
= T$, we denote by $\filt{s} \eqdef\post{s}{s}$ the filtering
distribution at time~$s$.

\section{Algorithms}
\label{sec:FFBS}
Conditionally on the observations $\chunk{y}{0}{T}$, the state
sequence $\{ X_s \}_{s \geq0}$ is a
time inhomogeneous Markov chain. This property remains true in the
\emph{time-reversed} direction.
Denote by $\BK{\eta}$ the so-called \emph{backward kernel} given by,
for any probability measure $\eta$ on $(\Xset, \Xsigma)$,
%
%
\begin{equation} \label{eq:backward-kernel}
\BK{\eta}(x, h) \eqdef\frac{\int\eta(\rmd x') \ensuremath{m}(x',
x) h(x')}
{\int\eta(\rmd x') \ensuremath{m}(x', x)} ,\qquad  h \in\mcb{\Xset
} .
\end{equation}
The posterior distribution $\post{s:T}{T}$
may be expressed as, for any integers $T > 0$, $s \in\{0, \dots, T-1\}
$ and any $h \in\mcb{\Xset^{T-s+1}}$,
%
%
\begin{equation} \label{eq:smoothing:backw_decomposition}
\quad \post{s:T}{T}(h) = \idotsint\filt{T}(\rmd x_T) \BK{\filt
{T-1}}(x_T, \rmd x_{T-1}) \cdots\BK{\filt{s}}(x_{s+1}, \rmd x_s)
h(\chunk{x}{s}{T}) .\hspace*{-10pt}
\end{equation}
Therefore, the joint smoothing distribution may be computed
recursively, backward in time, according to
%
%
\begin{equation} \label{eq:smoothing:backw_decomposition_recursion}
\post{s:T}{T}(h) = \idotsint\BK{\filt{s}}(x_{s+1}, \rmd x_s) \post
{s+1:T}{T}(\rmd\chunk{x}{s+1}{T}) h(\chunk{x}{s}{T}) .
\end{equation}

\subsection{The forward filtering backward smoothing algorithm}

As mentioned in the \hyperref[sec1]{Introduction}, the method proposed by \cite
{huerzelerkuensch1998,doucetgodsillandrieu2000} for approximating the
smoothing distribution is a two pass procedure.
In the forward pass, particle approximations $\filt[hat]{s}$ of the
filter distributions $\filt{s}$ are computed recursively for all time
steps from $s = 0$ up to $s = T$. The filter distribution flow $\{
\filt{s} \}_{s \geq0}$ satisfies the forward recursion
%
%
\begin{equation} \label{eq:forward-filtering-recursion}
\quad \filt{s}(h) = \frac{\unfilt{s}(h)}{\unfilt{s}(\1)}  \qquad \mbox{where } \unfilt{0}(h) = \Xinit(g_0 h) , \unfilt{s}(h) \eqdef
\unfilt{s-1} \M(g_s h) , s \geq1 ,
\end{equation}
for $h \in\mcb{\Xset}$, with $\1$ being the unity function $x
\mapsto1$ on $\Xset$. In terms of SMC, each filter distribution
$\filt{s}$ is approximated by means of a set of particles $\{ \epart
{s}{i} \}_{i = 1}^N$ and associated importance weights $\{ \ewght
{s}{i} \}_{i = 1}^N$ according to
%
%
\begin{equation} \label{eq:approximation-filtering-distribution}
\filt[hat]{s}(h) \eqdef\frac{\unfilt[hat]{s}(h)}{\unfilt[hat]{s}(\1
)} \qquad  \mbox{where } \unfilt[hat]{s}(h) \eqdef N^{-1} \sum_{i =
1}^N \ewght{s}{i} h(\epart{s}{i}) .
\end{equation}
%

Having produced, using methods described in Section \ref{section:APF}
below, a sequence of such weighted samples $\{ (\epart{t}{i}, \ewght
{t}{i}) \}_{i = 1}^N$, $1 \leq t \leq T$, an approximation of the
smoothing distribution is constructed in a backward pass by replacing,
in \eqref{eq:smoothing:backw_decomposition}, the filtering
distribution by its particle approximation. This yields
%
%
\begin{equation} \label{eq:smoothing:backw_decomposition_sample}
\post[hat]{s:T}{T}(h) \eqdef\int\!\!\cdots\!\!\int\!\filt[hat]{T}(\rmd x_T) \BK
{\filt[hat]{T-1}}(x_T,\rmd x_{T-1})\cdots\BK{\filt
[hat]{s}}(x_{s+1}, \rmd x_s) h(\chunk{x}{s}{T})\hspace*{-30pt}
\end{equation}
for any $h \in\mcb{\Xset^{T-s+1}}$. The approximation above can be
computed recursively in the backward direction according to
%
%
\begin{equation} \label{eq:smoothing:backw_decomposition_recursion_sample}
\post[hat]{s:T}{T}(h) = \int\!\!\cdots\!\!\int\BK{\filt[hat]{s}}(x_{s+1}, \rmd
x_s) \post[hat]{s+1:T}{T}(\rmd\chunk{x}{s+1}{T}) h(\chunk{x}{s}{T}).
\end{equation}
Now, by definition,
\[
\BK{\filt[hat]{s}}(x, h) = \sum_{i = 1}^N \frac{\ewght{s}{i}
\ensuremath{m}
(\epart{s}{i}, x)}{\sum_{\ell= 1}^N \ewght{s}{\ell} \ensuremath
{m}(\epart
{s}{\ell},x)} h ( \epart{s}{i} ) ,\qquad  h \in\mcb{\Xset} ,
\]
and inserting this expression into \eqref{eq:smoothing:backw_decomposition_sample} gives
%
%
\begin{equation} \label{eq:forward-filtering-backward-smoothing}
\post[hat]{s:T}{T}(h) \!=\! \sum_{i_s = 1}^N\!\!\cdots\!\!\sum_{i_T = 1}^N
\Biggl(\!\prod_{u=s+1}^T \frac{\ewght{u-1}{i_{u-1}} \ensuremath{m}(\epart
{u-1}{i_{u-1}},\epart{u}{i_u})}{\sum_{\ell=1}^N \ewght{u-1}{\ell}
\ensuremath{m}(\epart{u-1}{\ell},\epart{u}{i_u})} \!\Biggr) \frac{\ewght
{T}{i_T}}{\sumwght{T}} h (\epart{s}{i_s}, \dots, \epart{T}{i_T} )
,\hspace*{-30pt}
\end{equation}
of $\post{s:T}{T}(h)$, where $h \in\mcb{\Xset^{T-s-1}}$ and
%
%
\begin{equation} \label{eq:defOmega}
\sumwght{t} \eqdef\sum_{i=1}^N \ewght{t}{i} .
\end{equation}
The estimator $\post[hat]{s:T}{T}$ is impractical since the
cardinality of its support grows exponentially with the number $T-s$ of
time steps; nevertheless, it plays a~key role in the theoretical
developments that follow. A more practical approximation of this
quantity will be defined in the next section.
When the dimension of the input space is moderate, the computational
cost of evaluating the estimator can be reduced to $O(N \log N)$ by
using the \emph{fast
multipole method} as suggested in \cite
{klaasbriersdefreitasdoucetmaskelllang2006}; note, however, that this
method involves approximations that introduce some bias. On the other
hand, in certain specific scenarios, such as discrete Markov chains
with sparse transition matrices over large state spaces, the complexity
can even be reduced to $O(N T)$ without any truncation; see \cite
{barembruchgariviermoulines2009}.

\subsection{The forward filtering backward simulation algorithm}
\label{subsec:FFBSi}
The estimator~\eqref{eq:forward-filtering-backward-smoothing} may be
understood alternatively by noting that the normalized smoothing
weights define a probability distribution on the set $\{1, \dots, N\}
^{T - s}$ of trajectories associated with an inhomogeneous Markov
chain. Indeed, consider, for $t \in\{0, \dots, T-1\}$, the Markov
transition matrix $\{ \Lambda_t^N(i,j) \}_{i,j = 1}^N$ given by
%
%
\begin{equation} \label{eq:definition-transition-matrix-W}
\Lambda^N_{t}(i,j) = \frac{\ewght{t}{j} \ensuremath{m}(\epart
{t}{j}, \epart
{t+1}{i})}{\sum_{\ell=1}^N \ewght{t}{\ell} \ensuremath{m}(\epart
{t}{\ell},
\epart{t+1}{i})} ,\qquad  (i,j) \in\{1, \dots, N\}^2 .
\end{equation}
For $1 \leq t \leq T$, denote by
%
%
\begin{equation} \label{eq:definition-mcf}
\mcf{t}{N} \eqdef\sigma\{\chunk{Y}{0}{T}, (\epart{s}{i}, \ewght
{s}{i}); 0 \leq s \leq t, 1 \leq i \leq N\}
\end{equation}
the $\sigma$-algebra generated by the observations from time $0$ to
time $T$ as well as the particles and importance weights produced in
the forward pass up to time $t$.
The transition probabilities defined in \eqref
{eq:definition-transition-matrix-W} induce an inhomogeneous Markov
chain $\{ J_{u} \}_{u = 0}^T$ evolving backward in time as follows. At
time $T$, the random index $J_T$ is drawn from the set $\{1, \dots, N\}
$ such that $J_T$ takes the value $i$ with a probability proportional
to $\ewght{T}{i}$.
At time $t \leq T-1$ and given that the index $J_{t + 1}$ was drawn at
time step $t + 1$, the index $J_t$ is drawn from the set $\{1,\dots,N\}$
such that $J_t$ takes the value $j$ with probability $\Lambda
^N_{t}(J_t, j)$.
The joint distribution of $\chunk{J}{0}{T}$ is\vspace*{1pt} therefore given by, for
$\chunk{j}{0}{T} \in\{1, \dots, N\}^{T + 1}$,
%
%
\begin{equation} \label{eq:distribution-non-homogeneous}
\mathbb{P} [ \chunk{J}{0}{T} = \chunk{j}{0}{T} | \mcf{T}{N} ] =
\frac{\ewght{T}{j_T}}{\sumwght{T}} \Lambda_T^N(J_T, j_{T-1}) \cdots
\Lambda_0^N(j_1, j_0) .
\end{equation}
Thus, and this is a key observation, the FFBS estimator \eqref
{eq:forward-filtering-backward-smoothing} of the joint smoothing
distribution may be written as the conditional expectation
%
%
\begin{equation} \label{eq:EspCond}
\post[hat]{0:T}{T}(h) = \CPE{h (\epart{0}{J_0}, \dots, \epart
{T}{J_T} )}{\mcf{T}{N}} ,\qquad  h \in\mcb{\Xset^{T+1}} .
\end{equation}
We may therefore construct an unbiased estimator of the FFBS estimator
by drawing, conditionally independently given $\mcf{T}{N}$, $N$ paths
of $\{\chunk{J}{0}{T}^\ell\}_{\ell=1}^N$ of the inhomogeneous Markov
chain introduced above and then forming the (practical) estimator
%
%
\begin{equation} \label{eq:FFBSi:estimator}
\post[tilde]{0:T}{T}(h) = N^{-1} \sum_{\ell= 1}^N h ( \epart
{0}{J_0^\ell}, \dots, \epart{T}{J_T^\ell} ) ,\qquad  h \in\mcb
{\Xset^{T+1}} .
\end{equation}
This practical estimator was introduced in \cite
{godsilldoucetwest2004} (Algorithm 1, page 158).
For ease of notation, we have here simulated $N$ replicates of the
backward, index-valued Markov chain, but it would of course also be
possible to sample a~number of
paths that is either larger or smaller than $N$. The estimator\vspace*{-1pt} $\post
[hat]{0:T}{T}$ may be seen as a Rao--Blackwellized version of $\post
[tilde]{0:T}{T}$. The variance of the latter is increased, but the gain
in computational complexity is significant. The associated algorithm is
referred in the sequel to as the forward filtering backward simulation
(FFBSi) algorithm. In Section \ref{sec:TimeUniformExponentialFFBS},
forgetting properties of the inhomogeneous backward chain will play a
key role when establishing time uniform stability properties
of the proposed smoothing algorithm.

The computational complexity for sampling a single path of\vspace*{1pt} $\chunk
{J}{0}{T}$ is $O(NT)$; therefore, the overall computational effort\vspace*{-1pt}
spent when estimating $\post[tilde]{0:T}{T} $ using the FFBSi sampler
is $O(N^2T)$. Following \cite
{klaasbriersdefreitasdoucetmaskelllang2006}, this complexity can be
reduced further to $O(N \log(N) T)$ by means of the fast multipole
method; however, here again computational work is gained at the cost of
introducing additional approximations.

\subsection{A fast version of the forward filtering backward
simulation algorithm}
We are now ready to describe one of the main contributions of this
paper, namely a novel version of the FFBSi algorithm that can be proved
to reach linear computational complexity under appropriate assumptions.
At the end of the filtering phase of the FFBSi algorithm, all weighted
particle samples $\{(\epart{s}{i}, \ewght{s}{i})\}_{i=1}^N$, $0 \leq
s \leq T$, are available, and it remains to sample efficiently index
paths $\{ \chunk{J}{0}{T}^\ell\}_{\ell= 1}^N$ under the distribution
\eqref{eq:distribution-non-homogeneous}.
When the transition kernel $\ensuremath{m}$ is bounded from above in
the sense
that $\ensuremath{m}(x,x') \leq\sigma_+$ for all $(x,x') \in\Xset
\times\Xset
$, the paths can be simulated recursively backward in time using the
following accept--reject procedure. As in the standard FFBSi algorithm,
the recursion is initiated by sampling $J_T^1, \dots, J_T^N$
multinomially with probabilities proportional to $\{ \ewght{T}{i} \}
_{i = 1}^N$. For $s \in\{0, \dots, T\}$, let~$\mathcal{G}_s^N$ the
smallest $\sigma$-field containing $\mcf{T}{N}$ and $\sigma(J_t^\ell
\dvtx 1\leq l\leq N, t\geq s)$; then in order to draw $J_s^\ell$
conditionally on $\mathcal{G}_{s+1}^N$, we draw, first, an index
proposal~$I^{\ell}_s$ taking the value $i \in\{1, \ldots, N\}$ with
a probability proportional to~$\ewght{t}{i}$ and, second, an
independent uniform random variable $U_s^{\ell}$ on $[0,1]$. Then we\vspace*{-2pt}
set $J_s^\ell= I_s^\ell$ if $U_s^\ell\leq\ensuremath{m}(\epart
{s}{I_s^\ell},
\epart{s+1}{J_{s+1}^\ell}) / \sigma_+$; otherwise, we reject the
proposed index and make another trial.
To create samples of size $n \in\{1, \dots, N\}$ from a~multinomial
distribution on a set of $N$ elements at lines 
1 and 6, Algorithm~\ref{alg:smooth} relies on an efficient procedure
described in Appendix \ref{subsec:multisampling} that requires $O (
n(1+\log(1+N/n)) )$ elementary operations; see Proposition \ref
{prop:multisample}. Using this technique, the computational complexity
of Algorithm \ref{alg:smooth} can be upper-bounded as follows.

\begin{algorithm}[t]
\begin{algorithmic}[1]
\caption{FFBSi-smoothing}\label{alg:smooth}
\State\label{alg:samp1}sample $J_T^1, \dots, J_T^N$ multinomially
with probabilities proportional to $\{ \ewght{T}{i} \}_{i = 1}^N$
\For{$s$ from $T-1$ down to $0$}
\State$L \gets(1, \dots, N)$
\While{$L$ is not empty}
\State$n\gets$ size(L) \label{alg:updatem}
\State\label{alg:samp2} sample $I_1, \dots, I_n$ multinomially with
probabilities proportional to\mbox{}\hspace*{30pt}\qquad  $\{ \ewght{s}{i} \}_{i = 1}^N$
\State sample $U_1,\dots,U_n$ independently and uniformly over $[0,1]$
\State$nL \gets\varnothing$
\For{$k$ from $1$ to $n$}
\If{$U_k \leq \ensuremath{m}(\epart{s}{I(k)}, \epart
{s+1}{J_{s+1}^{L(k)}} )
/ \sigma_+$}
\State$J_{s}^{L(k)} \gets I_k$
\Else
\State$nL \gets nL \cup\{L(k)\}$
\EndIf
\EndFor
\State$L \gets nL$
\EndWhile
\EndFor
\end{algorithmic}
\end{algorithm}

For the bootstrap particle filter as well as the fully adapted
auxiliary particle filter (see Section \ref{section:APF} for precise
descriptions of these SMC filters), it is possible to derive an
asymptotic expression for the number of simulations required at line~8
of Algorithm \ref{alg:smooth} even if the kernel $m$ is not bounded
from below. The following result is obtained using theory derived in
the coming section.
\begin{prop} \label{prop:complexityBoostrap}
Assume that the transition kernel is bounded from above,
$m(x, x') \leq\sigma_+$ for all $(x, x') \in\Xset\times\Xset$.
At each iteration $s \in\{0, \dots, T-1\}$, let~$Z_s^N$ be the number
of simulations required in the accept--reject procedure of Algorithm~\ref{alg:smooth}.
\begin{itemize}
\item For the bootstrap auxiliary filter, $Z_s^N /N$ converges in
probability to
\[
\alpha(s) \eqdef\sigma_+ \post{s}{s-1}(g_s) \frac{\idotsint\rmd
x_{s+1} \prod_{u = s + 2}^T \int m(x_{u-1}, \rmd x_u)
g_u(x_u)}{\idotsint\post{s}{s-1}(\rmd x_s) g_s(x_s) \prod_{u = s +
1}^T m(x_{u-1}, \rmd x_u) g_u(x_u)}
\]
as $N$ goes to infinity.
\item In the fully adapted case, $Z_s^N / N$ converges in probability to
\[
\beta(s) \eqdef\sigma_+ \frac{\idotsint\rmd x_{s+1}\prod
_{u=s+2}^T \int m(x_{u-1}, \rmd x_u) g_u(x_u)}{\idotsint\filt{s}(\rmd
x_s) g_s(x_s) \prod_{u = s + 1}^T m(x_{u-1}, \rmd x_u) g_u(x_u)}
\]
as $N$ goes to infinity.
\end{itemize}
{\spaceskip=0.2em plus 0.05em minus 0.02em A sufficient condition for ensuring finiteness of $\alpha(s)$ and
$\beta(s)$ is that\break $\int g_u(x_u)\, \rmd x_u <$} $ \infty$ for all $u \geq0$.
\end{prop}
%

If the transition kernel satisfies stronger mixing conditions, it is
possible to derive an upper-bound on the computational complexity of
the FFBSi for any auxiliary particle filter, that is,\ the total number
of computations (and not only the total number of simulations). Note
that this result is not limited to the bootstrap and the fully adapted cases.


\begin{prop}\label{prop:complexity}
Assume that the transition kernel is bounded from below and above, that
is,\ $\sigma_- \leq m(x, x') \leq\sigma_+$ for all $(x, x') \in
\Xset\times\Xset$.
Let $C(N,T)$ denote the number of elementary operations required in
Algorithm \ref{alg:smooth}.
Then, there exists a constant $K$ such that such that
$ \PE[C(N,T)] \leq K N T \sigma_+ / \sigma_- $.
\end{prop}


The proofs of Propositions \ref{prop:complexityBoostrap} and \ref
{prop:complexity} involve theory developed in the coming section and
are postponed to Section \ref{sec:complexity:proofs}.

Before concluding this section on reduced complexity, let us mention
that efficient smoothing strategies have been considered by \cite
{fearnhead2005} using quasi-{M}onte {C}arlo methods. The smoother
(restricted to be one-dimensional) presented in this work has a
complexity that grows quadraticly in the number of particles $N$;
nevertheless, since the variance of the same decays as $O(N^{-2})$ (or
faster) thanks to the use of quasi-random numbers, the method is
equivalent to methods with complexity growing linearly in $N$ [since
the standard {M}onte {C}arlo variance is $O(N^{-1})$]. This solution is
of course attractive; we are however not aware of extensions of this
approach to multiple dimensions.

\subsection{Auxiliary particle filters}
\label{section:APF}
It remains to describe in detail how to produce sequentially the
weighted samples $\{(\epart{s}{i}, \ewght{s}{i})\}_{i = 1}^N$, $0
\leq s \leq T$, which can be done in several different ways (see \cite
{doucetdefreitasgordon2001,liu2001,cappemoulinesryden2005} and the
references therein). Still, most algorithms may be formulated within
the unifying framework of the {\em auxiliary particle filter} described
in the following. Let $\{\epart{0}{i}\}_{i = 1}^N$ be \iid\ random\vspace*{-3pt}
variables such that $\epart{0}{i} \sim\XinitIS{0}$ and set $\ewght
{0}{i} \eqdef\rmd\Xinit/ \rmd\XinitIS{0}(\epart{0}{i}) g_0(\epart
{0}{i})$. The weighted sample $\{ (\epart{0}{i}, \ewght{0}{i})\}_{i =
1}^N$ then targets the initial filter $\filt{0}$ in the sense that
$\filt[hat]{0}(h)$ estimates $\filt{0}(h)$ for $h \in\mcb{\Xset}$.
In order to describe the sequential structure of the auxiliary particle
filter, we proceed inductively and assume that we have at hand a
weighted sample $\{ (\epart{s-1}{i}, \ewght{s-1}{i}) \}_{i = 1}^N$
targeting $\filt{s - 1}$ in the same sense. Next, we aim at simulating
new particles from the target $\filt[tar]{s}$ defined as
%
%
\begin{equation} \label{eq:target-forward-filtering}
\filt[tar]{s}(h) = \frac{\unfilt[hat]{s-1} \M(g_s h)}{\unfilt
[hat]{s-1} \M(g_s)} ,\qquad  h \in\mcb{\Xset} ,
\end{equation}
in order to produce an updated particle sample approximating the
subsequent filter $\filt{s}$. Following \cite{pittshephard1999}, this
may be done by considering the \emph{auxiliary} target distribution
%
%
\begin{equation} \label{eq:target-auxiliary}
\filt[aux]{s}{} (i, h) \eqdef\frac{\ewght{s-1}{i} \M(\epart
{s-1}{i}, g_s h)}{\sum_{\ell= 1}^N \ewght{s-1}{\ell} \M(\epart
{s-1}{\ell}, g_s h)} ,\qquad  h \in\mcb{\Xset} ,
\end{equation}
on the product space $\{ 1, \dots, N \} \times\Xset$ equipped with
the product $\sigma$-algebra $\mathcal{P}(\{ 1, \dots, N \}) \otimes
\Xsigma$.
By construction, $\filt[tar]{s}$ is the marginal distribution of~$\phi
_s^{N, \mathrm{a}}$ with respect to the particle index. Therefore, we
may approximate the target distribution $\filt[tar]{s}$ on $(\Xset,
\Xsigma)$ by simulating from the auxiliary distribution and then
discarding the indices. More specifically, we first simulate pairs $\{
(I_s^i, \epart{s}{i}) \}_{i = 1}^N$ of indices and particles from the
instrumental distribution
%
%
\begin{equation} \label{eq:instrumental-distribution-filtering}
\instrpostaux{s}{s}(i, h) \propto\ewght{s-1}{i} \adjfunc
{s}{s}{\epart{s-1}{i}} \Kiss{s}{s}(\epart{s-1}{i},h) ,\qquad  h \in
\mcb{\Xset} ,
\end{equation}
on the product space $\{1, \dots, N\} \times\Xset$, where $\{
\adjfunc{s}{s}{\epart{s-1}{i}} \}_{i = 1}^N$ are so-called \emph
{adjustment multiplier weights} and $\Kiss{s}{s}$
is a Markovian \emph{proposal} transition kernel. In the sequel, we
assume for simplicity that $\Kiss{s}{s}(x, \cdot)$ has, for any $x
\in\Xset$, a density $\kiss{s}{s}(x, \cdot)$
with respect to the reference measure $\nu$. For each draw $(I_s^i,
\epart{s}{i})$, $i = 1, \dots, N$, we compute the importance weight
%
%
\begin{equation} \label{eq:weight-update-filtering}
\ewght{s}{i} \eqdef\frac{\ensuremath{m}(\epart{s-1}{I_s^i},\epart{s}{i})
g_s(\epart{s}{i})}{\adjfunc{s}{s}{\epart{s-1}{I_s^i}} \kiss
{s}{s}(\epart{s-1}{I_s^i},\epart{s}{i})} ,
\end{equation}
such that $\ewght{s}{i} \propto\rmd\phi_s^{N, \mathrm{a}} / \rmd
\instrpostaux{s}{s}(I_s^i, \epart{s}{i})$,
and associate it to the corresponding particle position $\epart
{s}{i}$. Finally, the indices $\{I_s^i \}_{i = 1}^N$ are discarded
whereupon $ \{ (\epart{s}{i}, \ewght{s}{i} ) \}_{i=1}^N$ is taken as
an approximation of $\filt{s}$. The simplest choice, yielding to the
so-called \emph{bootstrap particle filter algorithm} proposed by \cite
{gordonsalmondsmith1993}, consists of setting, for all $x \in\Xset$,
$\adjfunc{s}{s}{x} \equiv1$ and $\kiss{s}{s}(x, \cdot) \equiv
\ensuremath{m}
(x, \cdot)$. A more appealing---but often\vspace*{-1pt} computationally
costly---choice consists of using the adjustment weights
$\adjfunc{s}{s}{x} \equiv\adjfunc[fully]{s}{s}{x} \eqdef\int
\ensuremath{m}(x,
x') g_s(x') \,\rmd x'$, $x \in\Xset$, and the proposal transition density
\[
\kiss[fully]{s}{s}(x, x') \eqdef\frac{\ensuremath{m}(x,x')
g_s(x')}{\adjfunc
[fully]{s}{s}{x}} ,\qquad  (x, x') \in\Xset\times\Xset.
\]
In this case, the auxiliary particle filter is referred to as \emph
{fully adapted}. Other choices are discussed in \cite
{doucmoulinesolsson2008} and \cite{cornebisemoulinesolsson2008}.

\section{Convergence of the FFBS and FFBSi algorithms}

In this section, the convergence of the FFBS and FFBSi algorithms are
studied. For these two algorithms, nonasymptotic Hoeffding-type
deviation inequalities and CLTs are obtained. We also introduce a
decomposition, serving as a basis for most results obtained in this
paper, of the error $\post[hat]{0:T}{T} - \post{0:T}{T}$ and some
technical conditions under which the results are derived.

For any function $f\dvtx \Xset^d \to\rset$, we define by $\esssup{f}
\eqdef\sup_{x\in\Xset^d} |f(x)|$ and $\oscnorm{f} \eqdef\sup
_{(x,x') \in\Xset^d \times\Xset^d} |f(x)-f(x')|$ the supremum and
oscillator norms, respectively. Denote $\bar{\nset} \eqdef\nset\cup
\{\infty\}$ and consider the following assumptions where $T$ is the
time horizon which can be either a finite integer or infinity.
\begin{assum}\label{assum:bound-likelihood}
For all $0 \leq t \leq T$, $g_t(\cdot) > 0$ and $\sup_{0 \leq t
\leq T} \esssup{g_t} <\infty$.
\end{assum}

Define for $t \geq0$ the importance weight functions
%
%
\begin{equation}
\label{eq:definition-weightfunction-forward}
\quad \ewghtfunc{0}(x) \eqdef\frac{\rmd\Xinit}{\rmd\XinitIS{0}}(x)
g_0(x) \quad \mbox{and}\quad  \ewghtfunc{t}(x,x') \eqdef\frac{\ensuremath{m}(x,x')
g_t(x')}{\adjfunc{t}{t}{x} \kiss{t}{t}(x,x')} ,\qquad  t \geq1 .
\end{equation}

\begin{assum}
\label{assum:borne-FFBS}
\ $\sup_{1 \leq t \leq T} \esssup{\adjfunc{t}{t}{}} < \infty$ and
$\sup_{0 \leq t \leq T} \esssup{\ewghtfunc{t}} < \infty$.
\end{assum}

The latter assumption is rather mild; it holds in particular under
Assumption \ref{assum:bound-likelihood} for the bootstrap filter ($\kiss
{t}{t}=\ensuremath{m}$ and $\adjfunc{t}{t}{}\equiv1$) and is automatically
fulfilled in the fully adapted case ($\ewghtfunc{t} \equiv1$).

The coming proofs are based on a decomposition of the joint smoothing
distribution that we introduce below. For $0 \leq t < T$ and $h \in
\mcb{\Xset^{T+1}}$, define the kernel $L_{t,T}: \Xset^{t+1}\times
\Xsigma^{\otimes T+1} \to[0,1]$ by
%
%
\begin{equation} \label{eq:defLt}
L_{t,T}(\chunk{x}{0}{t},h) \eqdef\idotsint\Biggl(\prod_{u=t+1}^{T} \M
(x_{u-1},\rmd x_{u}) g_{u}(x_{u}) \Biggr) h(\chunk{x}{0}{T})
\end{equation}
and set $L_{T,T}(\chunk{x}{0}{T},h)\eqdef h(\chunk{x}{0}{T})$. By
construction, for every $t \in\{0,\dots,T\}$, the joint smoothing
distribution may be expressed as
%
%
\begin{equation} \label{eq:smooth:recursion}
\post{0:T}{T}(h)=\frac{\post{0:t}{t}[L_{t,T}(\cdot,h)]}{\post
{0:t}{t}[L_{t,T}(\cdot,\1)]} .
\end{equation}
This expression extends the classical forward--backward decomposition to
the joint smoothing distribution; here $L_{t,T}(\cdot, h)$ plays the
role of the so-called backward variable. This suggests to decompose the
error $\post[hat]{0:T}{T}(h) - \post{0:T}{T}(h)$ as the following
telescoping sum:
\begin{eqnarray}  \label{eq:decomp_Smooth}
\post[hat]{0:T}{T}(h) - \post{0:T}{T}(h) &=& \frac{\filt
[hat]{0}[L_{0,T}(\cdot,h)]}{\filt[hat]{0}[L_{0,T}(\cdot,\1)]} -
\frac{\filt{0}[L_{0,T}(\cdot,h)]}{\filt{0}[L_{0,T}(\cdot,\1)]}
\nonumber\\[-8pt]\\[-8pt]
&&{} + \sum_{t=1}^{T} \biggl\{ \frac{\post[hat]{0:t}{t}[L_{t,T}(\cdot
,h)]}{\post[hat]{0:t}{t}[L_{t,T}(\cdot,\1)]} - \frac{\post
[hat]{0:t-1}{t-1}[L_{t-1,T}(\cdot,h)]}{\post
[hat]{0:t-1}{t-1}[L_{t-1,T}(\cdot,\1)]} \biggr\} .\hspace*{-16pt}\nonumber
\end{eqnarray}
The first term on RHS of the decomposition above can be easily dealt
with since $\filt[hat]{0}$ is a weighted empirical distribution
associated to \iid\ random variables.

To cope with the terms in the sum of the RHS in \eqref
{eq:decomp_Smooth}, we introduce some kernels (depending on the {\em
past} particles) that stress the dependence \wrt\ the {\em current}
particules. More precisely, $\post[hat]{0:t}{t}[L_{t,T}(\cdot,h)]$ is
expressed as
%
%
\begin{equation} \label{eq:randomDenom}
\post[hat]{0:t}{t}[L_{t,T}(\cdot,h)] = \filt[hat]{t}[\F
^N_{t,T}(\cdot,h)] = \frac{\unfilt[hat]{t}[\F^N_{t,T}(\cdot
,h)]}{\unfilt[hat]{t}(\1)} ,
\end{equation}
where the random kernels $\F^N_{t,T}: \Xset\times\Xsigma^{\otimes
(T+1)}\to[0,1]$ are defined by: for all $0 < t \leq T$, and $x_t \in
\Xset$,
%
%
\begin{equation} \label{eq:definition-Ft}
\F^N_{t,T}(x_t,h) \eqdef\idotsint\BK{\filt[hat]{t-1}}(x_t, \rmd
x_{t-1}) \cdots\BK{\filt[hat]{0}}(x_{1}, \rmd x_{0}) L_{t,T}(\chunk
{x}{0}{t}, h) ,\hspace*{-25pt}
\end{equation}
and
%
%
\begin{equation} \label{eq:definition-Fs-As}
\F^N_{0,T}(x,h) \eqdef L_{0,T}(x,h) .
\end{equation}
We stress that the kernels $\F^N_{t,T}$ depend on the particles and\vspace*{-1pt}
weights $(\epart{s}{i}, \ewght{s}{i})_{i=1}^N$, $0 \leq s \leq t-1$,
through the particle approximations $\filt[hat]{t-1}, \ldots, \filt
[hat]{0}$ of the filter distributions. When proving the CLT for the
FFBS algorithm, it will be crucial to establish that for any $h \in
\mcb{\Xset^{T+1}}$, $\F^N_{t,T}(\cdot, h)$ converges (see Lemma~%
\ref{lem:limLG} below), as the number $N$ of particles tends to
infinity, to a deterministic function $\F_{t,T}(\cdot, h)$ given by
%
%
\begin{equation} \label{eq:definition-Ft-lim}
\F_{t,T}(x_t, h) \eqdef\idotsint\BK{\filt{t-1}}(x_t, \rmd x_{t-1})
\cdots\BK{\filt{0}}(x_{1}, \rmd x_{0}) L_{t,T}(\chunk
{x}{0}{t},h).\hspace*{-25pt}
\end{equation}
In the sequel, the case $h=\1$ will be of particular importance; in
that case, $L_{t,T}(\chunk{x}{0}{t},\1)$ does not depend on $\chunk
{x}{0}{t-1}$, yielding
%
%
\begin{equation}\label{eq:Lone}
\F^N_{t,T}(x_t, \1) = \F_{t,T}(x_t, \1) = L_{t,T}(\chunk{x}{0}{t}, \1)
\end{equation}
for all $x_{0:t} \in\Xset^{t + 1}$. Using these functions, the
difference appearing in the sum in \eqref{eq:decomp_Smooth} may then
be rewritten as
\begin{eqnarray}\label{eq:definition-A}
\nonumber
&&\frac{\post[hat]{0:t}{t}[L_{t,T}(\cdot,h)]}{\post
[hat]{0:t}{t}[L_{t,T}(\cdot,\1)]} - \frac{\post
[hat]{0:t-1}{t-1}[L_{t-1,T}(\cdot,h)]}{\post
[hat]{0:t-1}{t-1}[L_{t-1,T}(\cdot,\1)]} \hspace*{-25pt}\\
&&\qquad  = \frac{1}{\unfilt[hat]{t}[\F^N_{t,T}(\cdot,\1)]} \biggl(\unfilt
[hat]{t}[\F^N_{t,T}(\cdot,h)] - \frac{\filt[hat]{t-1}[\F
^N_{t-1,T}(\cdot,h)]}{\filt[hat]{t-1}[\F^N_{t-1,T}(\cdot,\1)]}
\unfilt[hat]{t}[\F^N_{t,T}(\cdot,\1)] \biggr)\hspace*{-25pt} \\
&&\qquad  = \frac{N^{-1} \sum_{\ell=1}^N \ewght{t}{\ell} G^N_{t,T}(\epart
{t}{\ell}, h)}{N^{-1} \sum_{\ell=1}^N \ewght{t}{\ell} \F
_{t,T}(\epart{t}{\ell}, \1)} ,\nonumber\hspace*{-25pt}
\end{eqnarray}
where the kernel $G^N_{t,T}: \Xset\times\Xsigma^{T+1} \to[0,1]$ is
defined by, for $x \in\Xset$,
%
%
\begin{equation} \label{eq:definition-G}
G^N_{t,T}(x,h) \eqdef\F^N_{t,T}(x,h) - \frac{\filt[hat]{t-1}[\F
^N_{t-1,T}(\cdot,h)]}{\filt[hat]{t-1}[\F^N_{t-1,T}(\cdot,\1)]} \F
^N_{t,T}(x,\1) .
\end{equation}
Similarly to $\F^N_{t,T}(\cdot, h)$, the functions $G^N_{t,T}(\cdot,
h)$ depend on the past particles; it will however be shown (see Lemma
\ref{lem:limLG} below) that $G^N_{t,T}(\cdot, h)$ converges to the
deterministic function given by, for $x \in\Xset$,
%
%
\begin{equation} \label{eq:definition-G-lim}
G_{t,T}(x, h) \eqdef\F_{t, T}\bigl(x, h - \post{0:T}{T}(h)\bigr) .
\end{equation}
The key property of this decomposition is stated in the following
lemma.

\begin{lem} \label{lem:GstTiszeromean}
Assume that Assumptions \ref{assum:bound-likelihood}--\ref{assum:borne-FFBS}
hold for some $T < \infty$. Then, for any $0 \leq t\leq T$, the
variables $\{\ewght{t}{\ell} G^N_{t,T}(\epart{t}{\ell},h)\}_{\ell=
1}^N$ are, conditionally on the $\sigma$-field $\mcf{t-1}{N}$, \iid\
with zero mean. Moreover, there exists a constant $C$ (that may depend
on $t$ and $T$) such that, for all $N \geq1$, $\ell\in\{1, \dots,
N\}$, and $h \in\mcb{\Xset^{T+1}}$,
\[
| \ewght{t}{\ell} G^N_{t,T}(\epart{t}{\ell}, h) | \leq\supnorm
{\ewght{t}{}} |G^N_{t,T}(\epart{t}{\ell}, h) | \leq C \oscnorm{h}
.
\]
\end{lem}

\begin{pf}
By construction, all pairs of particles and weights of the weighted
sample $\{(\epart{t}{\ell}, \ewght{t}{\ell})\}_{\ell=1}^N$ are
\iid\ conditionally on the $\sigma$-field $\mcf{t-1}{N}$. This
implies immediately that the variables $\{\ewght{t}{\ell}
G^N_{t,T}(\epart{t}{\ell},h)\}_{\ell= 1}^N$ are also \iid\
conditionally on the same $\sigma$-field $\mcf{t-1}{N}$. We now show
that $\mathbb{E}[ \ewght{t}{1} G^N_{t,T}(\epart{t}{1}, h) | \mcf
{t-1}{N}] = 0$. Using the definition of $G^N_{t,T}$ and the fact that
$\filt[hat]{t-1}[\F^N_{t-1,T}(\cdot, h)]$ and $\filt[hat]{t-1}[\F
^N_{t-1,T}(\cdot,\1)]$ are $\mcf{t-1}{N}$-measurable, we have
\begin{eqnarray*}
&&\CPE{\ewght{t}{1} G^N_{t,T}(\epart{t}{1},h)}{\mcf{t-1}{N}} \\
&&\qquad =\CPE{\ewght{t}{1} \F^N_{t,T}(x,h)}{\mcf{t-1}{N}}- \frac{\filt
[hat]{t-1}[\F^N_{t-1,T}(\cdot,h)]}{\filt[hat]{t-1}[\F
^N_{t-1,T}(\cdot,\1)]} \CPE{\ewght{t}{1} \F^N_{t,T}(x,\1)}{\mcf
{t-1}{N}} ,
\end{eqnarray*}
which is equal to zero provided that the relation
%
%
\begin{equation} \label{eq:espCondA}
\CPE{\ewght{t}{1} \F^N_{t,T}(\epart{t}{1},h)}{\mcf{t-1}{N}}=\frac
{\filt[hat]{t-1}[\F^N_{t-1,T}(\cdot,h)]}{\filt[hat]{t-1}(\adjfunc
{t}{t}{})}
\end{equation}
holds for any $h \in\mcb{\Xset}$. We now turn to the proof of \eqref
{eq:espCondA}. Note that for any $f \in\mcb{\Xset}$,
\begin{eqnarray}  \label{eq:technique}
\CPE{\ewght{t}{1} f(\epart{t}{1})}{\mcf{t-1}{N}} &=& \frac{\sum
_{\ell=1}^N \ewght{t-1}{\ell} \int\M(\epart{t-1}{\ell}, \rmd x)
g_t(x) f(x)}{\sum_{\ell=1}^N \ewght{t-1}{\ell} \adjfunc
{t}{t}{\epart{t-1}{\ell}}} \nonumber\\[-8pt]\\[-8pt]
&=& \frac{\filt[hat]{t-1}[\M(\cdot,g_t f)]}{\filt[hat]{t-1}(\adjfunc
{t}{t}{})} .\nonumber
\end{eqnarray}
It turns out that \eqref{eq:espCondA} is a consequence of \eqref
{eq:technique} with $f(\cdot)=\F_{t,T}^N(\cdot,h)$, but since $\F
^N_{t-1,T}(\cdot,h)$ is in general different from $\M(\cdot, g_t \F
_{t,T}^N(\cdot,h))$, we have to prove directly that
%
%
\begin{equation} \label{eq:cequonveut}
\filt[hat]{t-1} [ \F^N_{t-1,T}(\cdot,h) ] = \filt[hat]{t-1} [ \M
(\cdot, g_t \F_{t,T}^N(\cdot,h)) ] .
\end{equation}
Write
\begin{eqnarray}\label{eq:technique1}
&&\filt[hat]{t-1} [\M(\cdot,g_t \F_{t,T}^N(\cdot,h)) ] \nonumber \\
&&\qquad = \sumwght{t}^{-1} \sum_{\ell=1}^N \ewght{t-1}{\ell} \idotsint
\ensuremath{m}
(\epart{t-1}{\ell},x_t) g_t(x_t) \Biggl(\prod_{u=1}^{t} \BK{\filt
[hat]{u-1}}(x_u, \rmd x_{u-1}) \Biggr)\\
&&\hphantom{\qquad = \sumwght{t}^{-1} \sum_{\ell=1}^N \ewght{t-1}{\ell} \idotsint}
{}\times L_{t,T}(\chunk{x}{0}{t},h) \,\rmd
x_t.\nonumber
\end{eqnarray}
To simplify the expression in the RHS, we will use the two following equalities:
\begin{eqnarray}\label{eq:backRelat}
\Biggl(\!\sum_{\ell=1}^N\! \ewght{t-1}{\ell}
\ensuremath{m}(\epart{t-1}{\ell},x_t)\!\!\Biggr)\!
\BK{\filt[hat]{t-1}}(x_t,\rmd x_{t-1})\! &=&\!\sum_{\ell=1}^N\! \ewght
{t-1}{\ell} \ensuremath{m}(x_{t-1},x_t) \delta_{\epart{t-1}{\ell
}}(\rmd x_{t-1})  , \hspace*{-30pt}\\
\label{eq:recurL}
\int\M(x_{t-1}, \rmd x_t) g_t(x_t)L_{t,T}(\chunk{x}{0}{t}, h) &=&
L_{t-1,T}(\chunk{x}{0}{t-1}, h) .
\end{eqnarray}
The first relation is derived directly from the definition $\eqref
{eq:backward-kernel}$ of the backward kernel, the second is a recursive
expression of $L_{t,T}$ which is straightforward from the definition
\eqref{eq:defLt}. Now, \eqref{eq:backRelat} and \eqref{eq:recurL}
allow for writing
\begin{eqnarray*}
&& \sum_{\ell=1}^N \ewght{t-1}{\ell} \idotsint m(\epart{t-1}{\ell
},x_t) g_t(x_t) \prod_{u=1}^{t} \BK{\filt[hat]{u-1}}(x_u, \rmd
x_{u-1}) L_{t,T}(\chunk{x}{0}{t},h) \,\rmd x_t \\
&&\qquad = \sum_{\ell=1}^N \ewght{t-1}{\ell} \idotsint\M(x_{t-1}, \rmd
x_t) g_t(x_t) \delta_{\epart{t-1}{\ell}}(\rmd x_{t-1}) \\
&&\hphantom{\qquad = \sum_{\ell=1}^N \ewght{t-1}{\ell} \idotsint}
{}\times\prod
_{u=1}^{t-1} \BK{\filt[hat]{u-1}}(x_u, \rmd x_{u-1}) L_{t,T}(\chunk
{x}{0}{t},h) \\
&&\qquad = \sum_{\ell=1}^N \ewght{t-1}{\ell} \idotsint\delta_{\epart
{t-1}{\ell}}(\rmd x_{t-1}) \prod_{u=1}^{t-1} \BK{\filt
[hat]{u-1}}(x_u, \rmd x_{u-1}) L_{t-1,T}(\chunk{x}{0}{t-1},h) \\
&&\qquad = \sum_{\ell=1}^N \ewght{t-1}{\ell} \F^N_{t-1}(\epart{t-1}{\ell
},h).
\end{eqnarray*}
By plugging this expression into \eqref{eq:technique1}, we obtain
\eqref{eq:cequonveut} from which \eqref{eq:espCondA} follows via
\eqref{eq:technique}. Finally, $\mathbb{E}[\ewght{t}{1}
G^N_{t,T}(\epart{t}{1},h) | \mcf{t-1}{N} ] = 0$. It remains to check
that the random variable $\ewght{t}{1} G^N_{t,T}(\epart{t}{1}, h)$ is
bounded. But this is immediate since
\begin{eqnarray}  \label{eq:borne-ewghtfuncxG}
|\ewght{t}{1} G^N_{t,T}(\epart{t}{1}, h) | &=& \esssup{\ewghtfunc{t}}
\biggl|\F^N_{t,T}(\cdot, h) - \frac{\filt[hat]{t-1}[\F
^N_{t-1,T}(\cdot, h)]}{\filt[hat]{t-1}[\F^N_{t-1,T}(\cdot, \1)]} \F
_{t,T}(\cdot,\1)\biggr|_{\infty} \nonumber\hspace*{-30pt}\\[-8pt]\\[-8pt]
&\leq&2 \esssup{\ewghtfunc{t}} \esssup{\F^N_{t,T}(\cdot,\1
)}\oscnorm{h} \leq2 \esssup{\ewghtfunc{t}} \esssup{L_{t,T}(\cdot,
\1)} \oscnorm{h} .\hspace*{-30pt}\nonumber
\end{eqnarray}
\end{pf}

\subsection{Exponential deviation inequality}
\label{sec:exponentialFFBS}

We first establish a nonasymptotic deviation inequality.
Considering \eqref{eq:definition-A}, we are led to prove a Hoeffding
inequality for ratios. For this purpose, we use the following
elementary lemma which will play a key role in the sequel. The proof is
postponed to Appendix~\ref{sec:proof:lem:inegEssentielle}.

\begin{lem} \label{lem:inegessentielle}
Assume that $a_N$, $b_N$ and $b$ are random variables defined on the
same probability space such\vadjust{\goodbreak} that there exist positive constants $\beta
$, $B$, $C$ and $M$ satisfying:
\begin{enumerate}[(III)]
\item[(I)]\hypertarget{item:inegEssentielle-borne}
$|a_N/b_N|\leq M$, $\PP$-\as\ and $b \geq\beta$, $\PP$-\as,
\item[(II)]\hypertarget{item:inegEssentielle-expo1}
for all $\epsilon>0$ and all $N\geq1$, $\PP[|b_N-b|>\epsilon]\leq
B \rme^{-C N \epsilon^2}$,
\item[(III)]\hypertarget{item:inegEssentielle-expo2}
for all $\epsilon>0$ and all $N\geq1$, $\PP[ |a_N|>\epsilon]\leq
B \rme^{-C N (\epsilon/M )^2}$.
\end{enumerate}
Then
\[
\PP\biggl( \biggl| \frac{a_N}{b_N} \biggr| > \epsilon\biggr) \leq B \exp \biggl(-C N \biggl(\frac
{\epsilon\beta}{2M} \biggr)^2 \biggr) .
\]
\end{lem}

\begin{thm} \label{thm:Hoeffding-FFBS}
Assume that Assumptions \ref{assum:bound-likelihood}--\ref{assum:borne-FFBS}
hold for some $T < \infty$. Then, there exist constants $0 < B$ and $C
< \infty$ (depending on $T$) such that for all~$N$, $\epsilon> 0$,
and all measurable functions $h \in\mcb{\Xset^{T+1}}$,
%
%
\begin{equation} \label{eq:Hoeffding-1}
\PP[ | \post[hat]{0:T}{T}(h) - \post{0:T}{T}(h) | \geq\epsilon]
\leq B \rme^{-C N \epsilon^2 / \oscnorm[2]{h}} .
\end{equation}
In addition,
%
%
\begin{equation}\label{eq:LGN-unnormalised-L}
N^{-1} \sum_{\ell=1}^N \ewght{t}{\ell} \F_{t,T}(\epart{t}{\ell
},\1) \plim_{N \rightarrow\infty} \frac{\filt{t-1}[\F
_{t-1,T}(\cdot,\1)]}{\filt{t-1}(\adjfunc{t}{t}{})} .
\end{equation}
\end{thm}

\begin{rem}
As a by-product, Theorem \ref{thm:Hoeffding-FFBS} provides an
exponential inequality for the particle approximation of the filter.
For any $h\in\mcb{\Xset}$, define the function $\chunk{h}{0}{T}:
\Xset^{T+1} \to\rset$ by $\chunk{h}{0}{T}(\chunk{x}{0}{T})=
h(x_T)$. By construction, $\post{0:T}{T}(\chunk{h}{0}{T})= \filt
{T}(h)$ and
$\post[hat]{0:T}{T}(\chunk{h}{0}{T})= \filt[hat]{T}(h)$. With this
notation, equation~\eqref{eq:Hoeffding-1} may be rewritten as
\[
\PP[ | \filt[hat]{T}(h) - \filt{T}(h) | \geq\epsilon] \leq B \rme
^{-C N \epsilon^2 / \oscnorm[2]{h}} .
\]
An inequality of this form was first obtained by \cite
{delmoralmiclo2000} (see also \cite{delmoral2004},\break Chapter~7).
\end{rem}

\begin{pf}
We prove \eqref{eq:Hoeffding-1} by induction on $T$ using the
decomposition \eqref{eq:decomp_Smooth}. Assume that \eqref
{eq:Hoeffding-1} holds at time $T - 1$, for $\post[hat]{0:T
- 1}{T - 1}(h)$. Let $h \in\mcb{\Xset^{T+1}}$ and assume without
loss of generality that $\post{0:T}{T}(h) = 0$. Then \eqref
{eq:smooth:recursion} implies that $\filt{0}[L_{0,T}(\cdot, h)] = 0$
and the first term of the decomposition \eqref{eq:decomp_Smooth} thus becomes
%
%
\begin{equation}\label{eq:initialRecurInegExpo}
\frac{\filt[hat]{0}[L_{0,T}(\cdot,h)]}{\filt[hat]{0}[L_{0,T}(\cdot
,\1)]} = \frac{N^{-1}\sum_{i=0}^N \frac{\rmd\Xinit}{\rmd\XinitIS
{0}}(\epart{0}{i}) g_0(\epart{0}{i}) L_{0,T}(\epart{0}{i},h)}{N^{-1}
\sum_{\ell=0}^N \frac{\rmd\Xinit}{\rmd\XinitIS{0}}(\epart
{0}{\ell}) g_0(\epart{0}{\ell})L_{0,T}(\epart{0}{\ell},\1)} ,
\end{equation}
where $\{\epart{0}{i}\}_{i = 1}^N$ are \iid\ random variables with
distribution $\XinitIS{0}$. We obtain an exponential inequality for
\eqref{eq:initialRecurInegExpo} by applying Lemma \ref
{lem:inegessentielle} with
\[
\cases{
\displaystyle a_N = N^{-1} \sum_{i=0}^N \frac{\rmd\Xinit}{\rmd\XinitIS
{0}}(\epart{0}{i}) g_0(\epart{0}{i}) L_{0,T}(\epart{0}{i},h) ,
\cr
\displaystyle b_N = N^{-1} \sum_{i=0}^N \frac{\rmd\Xinit}{\rmd\XinitIS
{0}}(\epart{0}{i}) g_0(\epart{0}{i}) L_{0,T}(\epart{0}{i},\1)
, \cr
b = \beta= \Xinit[g_0(\cdot) L_{0,T}(\cdot,\1)] .
}
\]
Condition \hyperlink{item:inegEssentielle-borne}{(I)} is trivially satisfied
and conditions \hyperlink{item:inegEssentielle-expo1}{(II)} and
\hyperlink{item:inegEssentielle-expo2}{(III)} follow from the Hoeffding inequality for
\iid\ variables.

By \eqref{eq:decomp_Smooth} and \eqref{eq:definition-A}, it is now
enough to establish an exponential inequality for
%
%
\begin{equation} \label{eq:relation0}
\frac{\post[hat]{0:t}{t}[L_{t,T}(\cdot,h)]}{\post
[hat]{0:t}{t}[L_{t,T}(\cdot,\1)]} - \frac{\post
[hat]{0:t-1}{t-1}[L_{t-1,T}(\cdot,h)]}{\post
[hat]{0:t-1}{t-1}[L_{t-1,T}(\cdot,\1)]} = \frac{N^{-1} \sum_{\ell
=1}^N \ewght{t}{\ell} G^N_{t,T}(\epart{t}{\ell},h)}{N^{-1} \sum
_{\ell=1}^N \ewght{t}{\ell} \F_{t,T}(\epart{t}{\ell},\1)} ,\hspace*{-25pt}
\end{equation}
where $0 < t \leq T$. For that purpose, we use again Lemma \ref
{lem:inegessentielle} with
%
%
\begin{equation} \label{eq:definition-b-1}
\cases{
\displaystyle a_N = N^{-1}\sum_{\ell=1}^N \ewght{t}{\ell} G^N_{t,T}(\epart
{t}{\ell},h) ,\cr
\displaystyle b_N = N^{-1} \sum_{\ell=1}^N \ewght{t}{\ell} \F_{t,T}(\epart
{t}{\ell},\1),\cr
\displaystyle b = \beta= \frac{\filt{t-1}[\F_{t-1,T}(\cdot,\1)]}{\filt
{t-1}(\adjfunc{t}{t}{})} .
}
\end{equation}
By considering the LHS of \eqref{eq:relation0}, $ |{a_N}/{b_N} |\leq
2 \esssup{h}$, verifying condition \hyperlink{item:inegEssentielle-borne}{(I)}
in Lemma \ref{lem:inegessentielle}. By Lemma \ref
{lem:GstTiszeromean}, Hoeffding's inequality implies that there exist
constants $B$ and $C$ such that for all $N$, $\epsilon> 0$, and all
measurable function $h \in\mcb{\Xset^{T+1}}$,
\begin{eqnarray*}
&&\PP\Biggl[ \Biggl|N^{-1}\sum_{\ell=1}^N \ewght{t}{\ell} G^N_{t,T}(\epart
{t}{\ell},h) \Biggr| \geq\epsilon\Biggr]\\
&&\qquad = \PE\Biggl[\mathbb{P}\Biggl[ \Biggl|N^{-1}\sum_{\ell=1}^N \ewght{t}{\ell}
G^N_{t,T}(\epart{t}{\ell},h) \Biggr| \geq\epsilon\Big|\mcf{t-1}{N}\Biggr] \Biggr] \leq
B \rme^{-C N \epsilon^2/\oscnorm[2]{h}},
\end{eqnarray*}
verifying condition \hyperlink{item:inegEssentielle-expo2}{(III)} in Lemma \ref
{lem:inegessentielle}. It remains to verify condition \hyperlink{item:inegEssentielle-expo1}{(II)}. Since the pairs of particles and weights
of the weighted sample $\{(\epart{t}{\ell},\ewght{t}{\ell})\}_{\ell
=1}^N$ are \iid\ conditionally on $\mcf{t-1}{N}$, Hoeffding's
inequality implies that
%
%
\begin{equation}\label{eq:relat0}
\PP\Biggl[ \Biggl| b_N - \mathbb{E}\Biggl[N^{-1} \sum_{\ell=1}^N \ewght{t}{\ell}\F
_{t,T}(\epart{t}{\ell},\1)\Big|\mcf{t-1}{N}\Biggr] \Biggr| \geq\epsilon\Biggr] \leq B
\rme^{-CN \epsilon^2} .
\end{equation}
Moreover, by \eqref{eq:technique}, \eqref{eq:Lone}, and the
definition \eqref{eq:defLt}, we have
\begin{eqnarray}  \label{eq:relat2}
&&\mathbb{E}\Biggl[N^{-1} \sum_{\ell=1}^N \ewght{t}{\ell} \F_{t,T}(\epart
{t}{\ell}, \1)\Big|\mcf{t-1}{N}\Biggr] - b \nonumber\\[-8pt]\\[-8pt]
&&\qquad = \frac{ \filt[hat]{t-1}[\F_{t-1,T}(\cdot,\1)]}{\filt
[hat]{t-1}(\adjfunc{t}{t}{})} - \frac{ \filt{t-1}[\F_{t-1,T}(\cdot
,\1)]}{\filt{t-1}(\adjfunc{t}{t}{})} = \frac{\filt
[hat]{t-1}(H)}{\filt[hat]{t-1}(\adjfunc{t}{t}{})} ,\nonumber
\end{eqnarray}
with $H(\cdot) \eqdef\F_{t-1,T}(\cdot, \1) - \filt{t-1}[\F
_{t-1,T}(\cdot,\1)] \adjfunc{t}{t}{}(\cdot) / \filt{t-1}(\adjfunc
{t}{t}{})$. To obtain an exponential deviation inequality for \eqref
{eq:relat2}, we apply again Lemma \ref{lem:inegessentielle} with
\[
\cases{
a'_N = \filt[hat]{t-1}(H) , \cr
b'_N=\filt[hat]{t-1}(\adjfunc{t}{t}{}) , \cr
b' = \beta'= \filt[]{t-1}(\adjfunc{t}{t}{}) .
}
\]
By using the inequality
\begin{eqnarray*}
&&\F_{t-1,T}(x_{t - 1}, \1) \\
&&\qquad = \adjfunc{t}{t}{}(x_{t - 1}) \int\frac{\ensuremath{m}(x_{t - 1}, x_t)
g_{t}(x_t)}{\adjfunc{t}{t}{}(x_{t - 1})\kiss{t}{t}(x_{t - 1}, x_t) }
\kiss{t}{t}(x_{t - 1}, x_t) \F_{t,T}(x_t, \1) \,\rmd x_t \\
&&\qquad \leq\adjfunc{t}{t}{}(x_{t - 1}) \supnorm{\ewghtfunc{t}} \supnorm
{\F_{t,T}(\cdot,\1)} ,
\end{eqnarray*}
we obtain the bound $| \filt[hat]{t-1}(H)/\filt[hat]{t-1}(\adjfunc
{t}{t}{})|\leq2 \supnorm{\ewghtfunc{t}} \supnorm{\F_{t,T}(\cdot,\1
)}$ which verifies condition \hyperlink{item:inegEssentielle-borne}{(I)}.
Now, since $t-1 \leq T-1$ and $\filt{t-1}(H)=0$, the induction
assumption implies that conditions \hyperlink{item:inegEssentielle-expo1}{(II)}
and \hyperlink{item:inegEssentielle-expo2}{(III)} are satisfied for $|b'_N-b'|$
and~$|a'_N|$. Hence, Lemma \ref{lem:inegessentielle} shows that
%
%
\begin{equation}\label{eq:relat1}
\PP\Biggl[ \Biggl| \mathbb{E}\Biggl[ N^{-1} \sum_{\ell=1}^N \ewght{t}{\ell} \F
_{t,T}(\epart{t}{\ell},\1)\Big|\mcf{t-1}{N}\Biggr] - b \Biggr| > \epsilon\Biggr] \leq B
\rme^{-C N \epsilon^2} .
\end{equation}
Finally, \eqref{eq:relat0} and \eqref{eq:relat1} ensure that
condition \hyperlink{item:inegEssentielle-expo1}{(II)} in Lemma \ref
{lem:inegessentielle} is satisfied and an exponential deviation
inequality for \eqref{eq:relation0} follows. The proof of \eqref
{eq:Hoeffding-1} is complete. The last statement \eqref
{eq:LGN-unnormalised-L} of the theorem is a consequence of \eqref
{eq:relat0} and \eqref{eq:relat1}.
\end{pf}

The exponential inequality of Theorem \ref{thm:Hoeffding-FFBS} may be
more or less immediately extended to the FFBSi estimator.

\begin{cor} \label{cor:Hoeffding-FFBSi}
Under the assumptions of Theorem \ref{thm:Hoeffding-FFBS} there exist
constants $0 < B$ and $C < \infty$ (depending on $T$) such that for
all $N$, $\epsilon> 0$, and all measurable functions $h$,
%
%
\begin{equation} \label{eq:Hoeffding-2}
\PP[ | \post[tilde]{0:T}{T}(h) - \post{0:T}{T}(h) | \geq\epsilon
] \leq B \rme^{-C N \epsilon^2 / \oscnorm[2]{h}},
\end{equation}
where $\post[tilde]{0:T}{T}(h)$ is defined in \eqref{eq:FFBSi:estimator}.
\end{cor}

\begin{pf}
Using \eqref{eq:EspCond} and the definition of $\post
[tilde]{s:T}{T}(h)$, we may write
\begin{eqnarray*}
&&\post[tilde]{0:T}{T}(h) - \post[hat]{0:T}{T}(h) \\
&&\qquad = N^{-1} \sum_{\ell=1}^N \bigl[ h (\epart{0}{J_0^\ell},\dots, \epart
{T}{J_T^\ell} ) - \CPE{h (\epart{0}{J_0},\dots,\epart{T}{J_T}
)}{\mcf{T}{N}} \bigr] ,
\end{eqnarray*}
which implies \eqref{eq:Hoeffding-2} by the Hoeffding inequality and
\eqref{eq:Hoeffding-1}.
\end{pf}
%
\subsection{Asymptotic normality}
\label{sec:CLTFFBS}
We now extend the theoretical analysis of the forward-filtering
backward-smoothing estimator \eqref
{eq:smoothing:backw_decomposition_sample} to a CLT. Consider the
following mild assumption on the proposal distribution.
\begin{assum} \label{assum:bound-proposal-kernel}
$\supnorm{\ensuremath{m}} < \infty$ and $\sup_{0 \leq t \leq T}
\supnorm{\kiss
{t}{t}} < \infty$.
\end{assum}

CLTs for interacting particle models have been established in \cite
{delmoralmiclo2000,delmoral2004,doucmoulines2008}; the application to
these results to auxiliary particle filters is presented in~\cite
{johansendoucet2008} and \cite{doucmoulinesolsson2008}, Theorem 3.2.
Here, we base our proof on techniques developed in~\cite
{doucmoulines2008} (extending \cite{chopin2004} and \cite
{kuensch2005}). As noted in the previous section, it turns out crucial
that $G^N_{t,T}(\cdot, h)$ converges to a deterministic function as $N
\to\infty$. This convergence is stated in the following lemma.

\begin{lem} \label{lem:limLG}
Assume Assumptions \ref{assum:bound-likelihood}--\ref
{assum:bound-proposal-kernel}. Then, for any $h \in\mcb{\Xset}$ and
\mbox{$x \in\Xset$},
\begin{eqnarray*}
\lim_{N \to\infty} \F^N_{t,T}(x,h) &=& \F_{t,T}(x,h) ,\qquad  \PP
\mbox{-\as}, \\
\lim_{N \to\infty} G^N_{t,T} (x,h) &=& G_{t,T}(x,h) ,\qquad  \PP
\mbox{-\as},
\end{eqnarray*}
where $\F^N_{t,T}$, $\F_{t,T}$, $G^N_{t,T}$ and $G_{t,T}$ are defined
in \eqref{eq:definition-Ft}, \eqref{eq:definition-Ft-lim}, \eqref
{eq:definition-G} and \eqref{eq:definition-G-lim}. Moreover, there
exists a constant $C$ (that may depend on $t$ and $T$) such that for
all $N \geq1$, $\ell\in\{1, \dots, N\}$, and $h \in\mcb{\Xset}$,
\[
| \ewght{t}{\ell} G_{t,T}(\epart{t}{\ell},h) | \leq\supnorm
{\ewght{t}{}} |G_{t,T}(\epart{t}{\ell}, h) |\leq C \oscnorm{h}
,\qquad  \PP\mbox{-\as}
\]
\end{lem}

\begin{pf*}{Proof of Lemma \ref{lem:limLG}}
Let $h \in\mcb{\Xset}$ and $x_t \in\Xset$.
By plugging \eqref{eq:backward-kernel} with $\eta=\filt[hat]{t-1}$
into the definition \eqref{eq:definition-Ft} of $\F^N_{t,T}(x_t,h)$,
we obtain immediately
\begin{eqnarray*}
&&\F^N_{t,T}(x_t, h) \\
&&\qquad = \frac{\idotsint\filt[hat]{t-1}(\rmd x_{t-1}) \prod_{u = 0}^{t-2}
\BK{\filt[hat]{u}}(x_{u+1},\rmd x_{u}) \ensuremath{m}(x_{t-1}, x_t)
L_{t,T}(\chunk{x}{0}{t},h)}{\int\filt[hat]{t-1}(\rmd x_{t-1})
\ensuremath{m}
(x_{t-1}, x_{t})} \\
&&\qquad = \frac{\post[hat]{0:t-1}{t-1}[H([\cdot,x_t])]}{\filt
[hat]{t-1}[\ensuremath{m}
(\cdot, x_t)]}  \qquad \mbox{with } H(\chunk{x}{0}{t})\eqdef
\ensuremath{m}
(x_{t-1}, x_t) L_{t,T}(\chunk{x}{0}{t}, h) .
\end{eqnarray*}
The convergence of $\F^N_{t,T}(\cdot, h)$ follows from Theorem \ref
{thm:Hoeffding-FFBS}. The proof of the convergence of $G^N_{t,T}(\cdot
, h)$ follows the same lines. Finally, the final statement of the lemma
is derived from Lemma \ref{lem:GstTiszeromean} and the almost sure
convergence of~$G^N_{t,T}(\cdot, h)$ to $G_{t,T}(\cdot, h)$.
\end{pf*}

Now, we may state the CLT with an asymptotic variance given by a finite
sum of terms involving the limiting kernel $G_{t,T}$.
\begin{thm} \label{thm:FFBS-CLT}
Assume Assumptions \ref{assum:bound-likelihood}--\ref
{assum:bound-proposal-kernel}. Then, for any $h\in\mcb{\Xset^{T+1}}$,
%
%
\begin{equation} \label{eq:clt}
\sqrt{N} \bigl(\post[hat]{0:T}{T}(h) - \post{0:T}{T}(h) \bigr) \dlim\mathcal
{N} (0, \asymVar{0:T}{T}{h} )
\end{equation}
with
\begin{eqnarray}\label{eq:expression-covariance}
\asymVar{0:T}{T}{h} &\eqdef&\frac{\XinitIS{0}[\ewghtfunc{0}^2(\cdot
) G_{0,T}^2(\cdot,h)]}{\XinitIS{0}^2 [\ewghtfunc{0}(\cdot) \F_{0,
T}(\cdot, \1)]} + \sum_{t=1}^T \frac{\filt{t-1}[\upsilon
_{t,T}(\cdot,h)] \filt{t-1}(\adjfunc{t}{t}{})}{\filt{t-1}^2[\F
_{t-1,T}(\cdot,\1)]},  \\
\label{eq:definition-upsilon}
\upsilon_{t,T}(\cdot, h) &\eqdef&\adjfunc{t}{t}{\cdot} \int\Kiss
{t}{t}(\cdot, \rmd x) \ewghtfunc{t}^2(\cdot, x) G^2_{t,T}(x, h).
\end{eqnarray}
\end{thm}

\begin{pf}
Without loss of generality, we assume that $\post{0:T}{T}(h)=0$. We
show that $\sqrt{N} \post[hat]{0:T}{T}(h)$ may be expressed as
%
%
\begin{equation} \label{eq:decomp_Smooth_1}
\sqrt{N} \post[hat]{0:T}{T}(h)= \sum_{t=0}^T \frac
{V^N_{t,T}(h)}{W^N_{t,T}} ,
\end{equation}
where the sequence of random vectors $[V^N_{0,T}(h), \dots,
V^N_{T,T}(h)]$ is asymptotically normal and $[W^N_{0,T}, \dots,
W^N_{T,T}]$ converge in probability to a deterministic vector. The
proof of \eqref{eq:clt} then follows from Slutsky's lemma. Actually,
the decomposition \eqref{eq:decomp_Smooth_1} follows immediately from
the backward decomposition~\eqref{eq:decomp_Smooth} by setting, for $t
\in\{1, \dots, T\}$,
\begin{eqnarray*}
V_{0,T}^N(h) &\eqdef& N^{-1/2} \sum_{\ell=1}^N \frac{\rmd\Xinit
}{\rmd\XinitIS{0}}(\epart{0}{\ell})g_0(\epart{0}{\ell})
G_{0,T}(\epart{0}{\ell},h) , \\
V_{t,T}^N(h) &\eqdef& N^{-1/2}\sum_{\ell=1}^N \ewght{t}{\ell}
G^N_{t,T}(\epart{t}{\ell},h) , \\
W_{0,T}^N &\eqdef& N^{-1} \sum_{\ell=1}^N \frac{\rmd\Xinit}{\rmd
\XinitIS{0}}(\epart{0}{\ell})g_0(\epart{0}{\ell}) \F_{0,T}(\epart
{0}{\ell},\1) , \\
W_{t,T}^N &\eqdef& N^{-1} \sum_{\ell=1}^N \ewght{t}{\ell} \F
_{t,T}(\epart{t}{\ell},\1) .
\end{eqnarray*}
The convergence
\begin{eqnarray*}
W_{0,T}^N &\plim_{N \to\infty}& \Xinit[g_0(\cdot) \F_{0,T}(\cdot
,\1) ] , \\
W_{t,T}^N &\plim_{N \to\infty}& \frac{\filt{t-1}[\F_{t-1,T}(\cdot
,\1)]}{\filt{t-1}(\adjfunc{t}{t}{})}
\end{eqnarray*}
of $[W^N_{0,T}, \dots, W^N_{T,T}]$ to a deterministic vector is
established immediately using~\eqref{eq:LGN-unnormalised-L} and noting
that the initial particles $(\epart{0}{i})_{i=1}^N$ are \iid\ We devote
the rest of the proof to showing that the sequence of random vectors
$[V^N_{0,T}(h), \dots,\break V^N_{T,T}(h)]$ is asymptotically normal.
Proceeding recursively in time, we prove by induction over $t \in
\{0, \dots, T\}$ (starting with $t=0$) that $[V^N_{0,T}(h), \dots,
V^N_{t,T}(h)]$ is asymptotically normal. More precisely, using the
Cram\'er--Wold device, it is enough to show that for all scalars
$(\alpha_0,\dots,\alpha_t) \in\rset^{t+1}$,
%
%
\begin{equation} \label{eq:multivariate-CLT-FFBS}
\sum_{r=0}^t \alpha_r V_{r,T}^N(h) \dlim_{N \to\infty} \mathcal
{N} \Biggl(0, \sum_{r=0}^t \alpha_r^2 \incrasymVar{r}{T}{h} \Biggr) ,
\end{equation}
where, for $r \geq1$,
\[
\incrasymVar{0}{T}{h}\eqdef\XinitIS{0}[\ewghtfunc{0}^2
G_{0,T}^2(\cdot,h)] ,\qquad  \incrasymVar{t}{T}{h} \eqdef\frac{\filt
{t-1}[\upsilon_{t,T}(\cdot,h)]}{\filt{t-1}(\adjfunc{t}{t}{})} .
\]
The case $t = 0$ is elementary since the initial particles $\{ \epart
{0}{i} \}_{i=1}^N$ are \iid\ Assume now that \eqref
{eq:multivariate-CLT-FFBS} holds for some $t-1 \leq T$; for
all scalars $(\alpha_1,\dots,\break\alpha_{t-1}) \in\rset^{t-1}$,
%
%
\begin{equation} \label{eq:induction-assumption}
\sum_{r=s}^{t-1} \alpha_r V_{r,T}^N(h) \dlim_{N \to\infty}
\mathcal{N} \Biggl(0, \sum_{r=s}^{t-1} \alpha_r^2 \incrasymVar{r}{T}{h} \Biggr)
.
\end{equation}
The sequence of random variable $V_{t,T}^N(h)$ may be expressed as an
additive function of a triangular array of random variables,
\[
V_{t,T}^N(h)= \sum_{\ell=1}^N U_{N,\ell} ,\qquad  U_{N,\ell}\eqdef
\ewght{t}{\ell} G^N_{t,T}(\epart{t}{\ell},h) / \sqrt{N} ,
\]
where $G^N_{t,T}(x,h)$ is defined in \eqref{eq:definition-G}. Lemma
\ref{lem:GstTiszeromean} implies that $\mathbb{E}[ V_{t,T}^N(h) |
\mcf{t-1}{N}] = 0$, yielding
\[
\mathbb{E}\Biggl[\sum_{r=0}^t \alpha_r V_{r,T}^N(h)\Big|\mcf{t-1}{N}\Biggr] = \sum
_{r=0}^{t-1} \alpha_r V_{r,T}^N(h) \dlim_{N \to\infty}
\mathcal{N} \Biggl(0, \sum_{r=1}^{t-1} \alpha_r^2 \incrasymVar{r}{T}{h} \Biggr),
\]
where the last limit follows by the induction assumption hypothesis
\eqref{eq:induction-assumption}. By \cite{doucmoulines2008}, Theorem A.3, page 2360, as the random variables $\{ U_{N,\ell} \}
_{\ell=1}^N$ are centered and conditionally independent given $\mcf
{t-1}{N}$, \eqref{eq:multivariate-CLT-FFBS}\vadjust{\goodbreak} holds provided that the
asymptotic smallness condition
%
%
\begin{equation} \label{eq:triangular-3}
\sum_{\ell=1}^N \mathbb{E}\bigl[U_{N,\ell}^2 \mathbh{1}_{\{|U_{N,\ell}|
\geq\epsilon\}}|\mcf{t-1}{N}\bigr] \plim_{N \to\infty} 0
\end{equation}
holds for any $\epsilon> 0$ and that the conditional variance converges:
%
%
\begin{equation} \label{eq:triangular-2}
\sum_{\ell=1}^N \CPE{U_{N,\ell}^2}{\mcf{t-1}{N}} \plim_{N \to
\infty} \incrasymVar{t}{T}{h} .
\end{equation}
Lemma \ref{lem:GstTiszeromean} implies that $|U_{N,\ell}|\leq C
\oscnorm{h} / \sqrt{N}$, verifying immediately the asymptotic
smallness condition \eqref{eq:triangular-3}. To conclude the proof, we
thus only need to establish the convergence \eqref{eq:triangular-2} of
the asymptotic variance. Via Lem\-ma~\ref{lem:GstTiszeromean} and
straightforward computations, we conclude that
\begin{eqnarray}  \label{eq:espCondU}
\sum_{\ell=1}^N \CPE{U^2_{N,\ell}}{\mcf{t-1}{N}} &=& \CPE{ ( \ewght
{t}{1} G^N_{t,T}(\epart{t}{1},h) )^2}{\mcf{t-1}{N}}\nonumber
\\
&=& \int\sum_{\ell=1}^N \frac{\ewght{t-1}{\ell} \adjfunc
{t}{t}{\epart{t-1}{\ell}} \Kiss{t}{t}(\epart{t-1}{\ell}, \rmd x)}
{\sum_{j=1}^N \ewght{t-1}{j} \adjfunc{t}{t}{\epart{t-1}{j}}} (
\ewghtfunc{t}(\epart{t-1}{\ell},x) G^N_{t,T}(x,h) )^2\hspace*{-20pt}
\nonumber\\[-8pt]\\[-8pt]
&=& \biggl(\frac{\sumwght{t-1}}{\sum_{j=1}^N \ewght{t-1}{j} \adjfunc
{t}{t}{\epart{t-1}{j}}} \biggr) \Biggl(\frac{1}{\sumwght{t-1}}\sum_{\ell=1}^N
\ewght{t-1}{\ell} \upsilon^N_{t,T}(\epart{t-1}{\ell},h) \Biggr) \nonumber\\
&=& \frac
{\filt[hat]{t-1}[\upsilon^N_{t,T}(\cdot,h)]}{\filt
[hat]{t-1}(\adjfunc{t}{t}{})} ,\nonumber
\end{eqnarray}
where $\sumwght{t}$ is defined in \eqref{eq:defOmega} and
\[
\upsilon^N_{t,T}(\cdot,h) \eqdef\adjfunc{t}{t}{\cdot} \int\Kiss
{t}{t}(\cdot, \rmd x) \ewghtfunc{t}^2(\cdot,x) [ G^N_{t,T}(x,h) ]^2
.
\]
The denominator in on RHS of \eqref{eq:espCondU} converges evidently
in probability to~$\filt{t-1}(\adjfunc{t}{t}{})$ by Theorem \ref
{thm:Hoeffding-FFBS}. The numerator is more complex since $\upsilon
^N_{t,T}$ depends on $G^N_{t,T}$ whose definition involves all the
approximations $\filt[hat]{t-1},\ldots,\filt[hat]{0}$ of the past\vspace*{-3pt}
filters. To obtain its convergence, note that, by Theorem \ref
{thm:Hoeffding-FFBS}, $\filt[hat]{t-1}(\upsilon_{t,T}(\cdot,h))
\plim\filt{t-1}(\upsilon_{t,T}(\cdot,h))$ as $N$ tends to infinity;
hence, it only remains to prove that
%
%
\begin{equation} \label{eq:technicos}
\filt[hat]{t-1} [\upsilon^N_{t,T}(\cdot, h) - \upsilon_{t,T}(\cdot
, h) ] \plim_{N \to\infty} 0 .
\end{equation}
For that purpose, introduce the following notation: for all $x \in
\Xset$,
\begin{eqnarray*}
A_N(x) &\eqdef&\filt[hat]{t-1} [\adjfunc{t}{t}{\cdot} \kiss
{t}{t}(\cdot, x) \ewghtfunc{t}^2(\cdot,x)| (G^N_{t,T}(x,h))^2 -
G_{t,T}^2(x,h)| ] , \\
B_N(x) &\eqdef&\filt[hat]{t-1} [\adjfunc{t}{t}{\cdot} \kiss
{t}{t}(\cdot, x) ] .
\end{eqnarray*}
Applying Fubini's theorem,
%
%
\begin{equation}\label{eq:intAn}
\lim_{N \to\infty} \PE\biggl[ \int A_N(x) \,\rmd x \biggr] = \lim_{N \to\infty
} \int\PE[A_N(x)] \,\rmd x = 0 ,
\end{equation}
where the last equality is due to the generalized Lebesgue convergence
theorem \cite{royden1988}, Proposition 18, page 270, with $f_N(x) = \PE
[A_N(x)]$ and $g_N(x) = 2 C \oscnorm{h} \PE[B_N(x)]$ provided that
the following conditions hold:
\begin{enumerate}[(iii)]
\item[(i)]\hypertarget{item:borneA} for any $x\in\Xset$, $\PE[A_N(x)] \leq2
C^2 \oscnorm[2]{h} \PE[B_N(x)]$,
\item[(ii)]\hypertarget{item:limA} for any $x \in\Xset$, $\lim_{N \to\infty}
\PE[A_N(x)] = 0$, $\PP$-\as,
\item[(iii)]\hypertarget{item:limIntB} $\lim_{N \to\infty} \int\PE[B_N(x)]
\,\rmd x=\int\lim_{N \to\infty} \PE[B_N(x)] \,\rmd x$.
\end{enumerate}

\textit{Proof of} \hyperlink{item:borneA}{(i)}. The bound follows directly from
Lemmas \ref{lem:limLG} and \ref{lem:GstTiszeromean}.

\textit{Proof of} \hyperlink{item:limA}{(ii)}. Using again Lemmas \ref{lem:limLG}
and \ref{lem:GstTiszeromean}, for any $x \in\Xset$,
\begin{eqnarray*}
A_N(x) &\leq&2C^2 \supnorm{\adjfunc{t}{t}{}} \supnorm{\kiss{t}{t}}
\oscnorm[2]{h} , \\
\limsup_{N \to\infty} A_N(x) &\leq&\supnorm{\adjfunc{t}{t}{}\kiss
{t}{t}\ewghtfunc{t}^2} \limsup_{N \to\infty} | (G^N_{t,T}(x, h))^2
- G_{t,T}^2(x,h) | = 0 ,\qquad  \PP\mbox{-\as}
\end{eqnarray*}
These two inequalities combined with $A_N(x)\geq0$ allow for applying
the Lebesgue dominated convergence theorem, verifying condition \hyperlink{item:limA}{(ii)}.

\textit{Proof of} \hyperlink{item:limIntB}{(iii)}. We have
%
\begin{eqnarray*}
\lim_{N \to\infty} \int\PE[B_N(x)] \,\rmd x &\stackrel{\mathrm{(a)}}{=} &\lim
_{N \to\infty} \PE\biggl[ \filt[hat]{t-1} \biggl( \adjfunc{t}{t}{\cdot} \int
\kiss{t}{t}(\cdot, x) \,\rmd x \biggr) \biggr] \\
&\stackrel{\mathrm{(b)}}{=}& \filt{t-1} (\adjfunc{t}{t}{} ) \stackrel{\mathrm{(c)}}{=}
\int\filt{t-1} (\adjfunc{t}{t}{\cdot} \kiss{t}{t}(\cdot, x) )
\,\rmd x \\
&\stackrel{\mathrm{(d)}}{=}& \int\lim_{N \to\infty} \PE[\filt
[hat]{t-1} (\adjfunc{t}{t}{\cdot}\kiss{t}{t}(\cdot, x) ) ] \,\rmd x\\
& =&
\int\lim_{N \to\infty} \PE[B_N(x)] \,\rmd x ,
\end{eqnarray*}
where (a) and (c) are consequences of Fubini's theorem and (b) and (d)
follows from the $\lone$-convergence of $\filt[hat]{t}(h)$ to $\filt
{t}(h)$ (see Theorem \ref{thm:Hoeffding-FFBS}) with $h(\cdot
)=\adjfunc{t}{t}{\cdot}$ and $h(\cdot)=\adjfunc{t}{t}{\cdot}\kiss
{t}{t}(\cdot, x)$.

Thus, \eqref{eq:intAn} holds, yielding that $\int A_N(x) \,\rmd x \plim
0$ as $N$ tends to infinity. This in turn implies \eqref{eq:technicos}
via the inequality
\[
| \filt[hat]{t-1} [\upsilon^N_{t,T}(\cdot,h)-\upsilon_{t,T}(\cdot
,h) ] | \leq\int A_N(x) \,\rmd x .
\]
This establishes \eqref{eq:multivariate-CLT-FFBS} and therefore
completes the proof.
\end{pf}

The weak convergence of $\sqrt{N} ( \post[hat]{0:T}{T}(h) - \post
{0:T}{T}(h))$ for the FFBS algorithm implies more or less immediately
the one of $\sqrt{N} (\post[tilde]{0:T}{T}(h) - \post{0:T}{T}(h) )$
for the FFBSi algorithm.
\begin{cor}
Under the assumptions of Theorem \ref{thm:FFBS-CLT},
\begin{eqnarray} \label{eq:cltSi}
&&\sqrt{N} \bigl(\post[tilde]{0:T}{T}(h) - \post{0:T}{T}(h) \bigr)
\nonumber\\[-8pt]\\[-8pt]
&&\qquad \dlim\mathcal{N} \bigl(0, \post{0:T}{T}^2 [h - \post{0:T}{T}(h) ] +
\asymVar{0:T}{T}{h - \post{0:T}{T}(h)} \bigr) .\nonumber
\end{eqnarray}
\end{cor}

\begin{pf}
Using \eqref{eq:EspCond} and the definition of $\post
[tilde]{0:T}{T}(h)$, we may write
\begin{eqnarray*}
&&\sqrt{N} \bigl(\post[tilde]{0:T}{T}(h) - \post{0:T}{T}(h) \bigr) \\
&&\qquad = N^{-1/2} \sum_{\ell=1}^N \bigl[ h (\epart{0}{J_0^\ell}, \dots, \epart
{T}{J_T^\ell} ) - \CPE{h (\epart{0}{J_0}, \dots, \epart{T}{J_T}
)}{\mcf{T}{N}} \bigr] \\
&&\qquad\quad {} + \sqrt{N} \bigl( \post[hat]{0:T}{T}(h) - \post{0:T}{T}(h) \bigr) .
\end{eqnarray*}
Note that since $\{ \chunk{J}{0}{T}^\ell\}_{\ell= 1}^N$ are \iid\
conditional on $\mcf{T}{N}$, \eqref{eq:cltSi} follows from \eqref
{eq:clt} and direct application of \cite{doucmoulines2008}, Theorem A.3, page 2360, by noting that
\begin{eqnarray*}
&&N^{-1} \sum_{\ell=1}^N \mathbb{E}\bigl[ \{h (\epart{0}{J_0^\ell},\dots,
\epart{T}{J_T^\ell} ) - \CPE{h (\epart{0}{J_0},\dots,\epart
{T}{J_T} )}{\mcf{T}{N}} \}^2|\mcf{T}{N}\bigr] \\
&&\qquad = \bigl(\post[hat]{0:T}{T} [h - \post[hat]{0:T}{T}(h) ] \bigr)^2
\plim\bigl(\post{0:T}{T} [h - \post{0:T}{T}(h) ] \bigr)^2 .
\end{eqnarray*}
\upqed
\end{pf}

\section{Time uniform bounds}
\label{sec:TimeUniformExponentialFFBS}
Most often, it is not required to compute the joint smoothing
distribution but rather the marginal smoothing distributions~$\post{s}{T}$.
Considering \eqref{eq:forward-filtering-backward-smoothing} for a
function $h$ that depends on the component~$x_s$ only, we obtain
particle approximations of the marginal smoothing distributions by
associating the set $\{ \epart{s}{j} \}_{j = 1}^N$ of particles with
weights obtained by marginalizing the joint smoothing weights according to
\[
\ewght{s|T}{i_s} = \sum_{i_{s+1} = 1}^N \cdots\sum_{i_T = 1}^N \prod
_{u=s+1}^t \frac{\ewght{u-1}{i_{u-1}} \ensuremath{m}(\epart{u-1}{i_{u-1}},
\epart{u}{i_u})}{\sum_{\ell=1}^N \ewght{u-1}{\ell} \ensuremath
{m}(\epart
{u-1}{\ell}, \epart{u}{i_u})} \frac{\ewght{T}{i_T}}{\sumwght{T}}.
\]
It is easily seen that these marginal weights may be recursively
updated backward in time as
%
%
\begin{equation}
\label{eq:FFBS-marginal-weight-update}
\ewght{s|T}{i} = \sum_{j = 1}^N \frac{
\ewght{s}{i} \ensuremath{m}(\epart{s}{i}, \epart{s+1}{j})}{\sum
_{\ell=1}^N
\ewght{s}{\ell} \ensuremath{m}(\epart{s}{\ell}, \epart{s+1}{j})}
\ewght
{s+1|T}{j} .
\end{equation}
In this section, we study the long-term behavior of the marginal
fixed-interval smoothing distribution estimator.
For that purpose, it is required to impose a type of mixing condition
on the Markov transition kernel; see~\cite
{chiganskylipstervanhandel2008} and the references therein. For
simplicity, we consider elementary but strong conditions which are
similar to the ones used in \cite{delmoral2004}, Chapter 7.4, or \cite
{cappemoulinesryden2005}, Chapter 4; these conditions, which points to
applications where the state space $\Xset$ is compact, can be relaxed,
but at the expense of many technical difficulties \cite
{chiganskylipster2004,vanhandel2008a,vanhandel2008b,vanhandel2009a}.

\begin{assum}
\label{assum:strong-mixing-condition}
There exist two constants $0 < \sigma_- \leq\sigma_+ < \infty$,
such that, for any $(x,x') \in\Xset\times\Xset$,
%
%
\begin{equation}
\label{eq:minorq}
\sigma_- \leq\ensuremath{m}(x,x') \leq\sigma_+.
\end{equation}
In addition, there exists a constant $c_- > 0$ such that, $\int\Xinit
(\rmd x_0) g_0(x_0) \geq c_-$ and for all $t \geq1$,
%
%
\begin{equation}
\label{eq:minorg}
\inf_{x \in\Xset} \int\M(x, \rmd x') g_t(x') \geq c_- > 0 .
\end{equation}
\end{assum}

Assumption \ref{assum:strong-mixing-condition} implies that $\nu
(\Xset)<\infty$; in the sequel, we will consider without loss of
generality that $\nu(\Xset)=1$.
Note also that, under Assumption \ref{assum:strong-mixing-condition}, the
average number of simulations required in the accept--reject mechanism
per sample of the FFBSi algorithm is bounded by $\sigma_+/\sigma_-$.

The goal of this section consists in establishing, under the
assumptions mentioned above, that the FFBS approximation of the \emph
{marginal} fixed interval smoothing probability satisfies an
exponential deviation inequality with constants that are uniform in
time and, under the same assumptions, that the variance of the CLT is
uniformly bounded in time.

For obtaining these results, we will need upper-bounds on $G^N_{t,T}$
and $G_{t,T}$ that are more precise than the ones stated in Lemmas \ref
{lem:GstTiszeromean} and \ref{lem:limLG}. For any function $h \in\mcb
{\Xset}$ and $s\leq T$, define the extension $\extens{h}{s}{T} \in
\mcb{\Xset^{T+1}}$ of $h$ to $\Xset^{T + 1}$ by
%
%
\begin{equation} \label{eq:defH}
\extens{h}{s}{T}(\chunk{x}{0}{T}) \eqdef h(x_s) ,\qquad  x_{0:T} \in
\Xset^{T + 1} .
\end{equation}

\begin{lem} \label{lem:G-uniform}
Assume that Assumptions \ref{assum:bound-likelihood}--\ref
{assum:strong-mixing-condition} hold with $T = \infty$. Let $s \leq
T$. Then, for all $t$,$T$, $N \geq1$, and $h \in\mcb{\Xset}$,
%
%
\begin{equation}\label{eq:time-unif-G}
\supnorm{G^N_{t,T}(\cdot, \extens{h}{s}{T})} \leq\rho
^{|t-s|}\oscnorm{h} \supnorm{\F_{t,T}(\cdot,\1)},
\end{equation}
where $\F_{t,T}$ is defined in \eqref{eq:definition-Ft-lim} and
%
%
\begin{equation}\label{eq:definition-rho}
\rho= 1 - \frac{\sigma_-}{\sigma_+} .
\end{equation}
Moreover, for all $t$, $T \geq1$, and $h \in\mcb{\Xset}$,
%
%
\begin{equation}\label{eq:time-unif-G-lim}
\supnorm{G_{t,T}(\cdot, \extens{h}{s}{T})} \leq\rho
^{|t-s|}\oscnorm{h} \supnorm{\F_{t,T}(\cdot,\1)}.
\end{equation}
\end{lem}

\begin{pf}
Using \eqref{eq:Lone} and \eqref{eq:definition-G},
%
%
\begin{equation}\label{eq:expressionG}
\frac{G^N_{t,T}(x, \extens{h}{s}{T})}{\F_{t,T}(x,\1)}=\frac{\F
^N_{t,T}(x,\extens{h}{s}{T})}{\F^N_{t,T}(x,\1)}- \frac{\filt
[hat]{t-1}[\F^N_{t-1,T}(\cdot,\extens{h}{s}{T})]}{\filt
[hat]{t-1}[\F^N_{t-1,T}(\cdot,\1)]} .
\end{equation}
To prove \eqref{eq:time-unif-G}, we will rewrite \eqref
{eq:expressionG} and obtain an exponential bound by either using
ergodicity properties of the ``a posteriori'' chain (when $t \leq s$),
or by using ergodicity properties of the backward kernel (when $t>s)$.

Assume first that $t \leq s$. The quantity $L_{t,T}(\chunk{x}{0}{t},
\extens{h}{s}{T})$ does not depend on $\chunk{x}{0}{t-1}$ so that by
\eqref{eq:definition-Ft} and definition \eqref{eq:defLt} of $L_{t,T}$,
\begin{eqnarray}\label{eq:one}
\F^N_{t,T}(x_t,\extens{h}{s}{T}) &=& L_{t,T}(\chunk{x}{0}{t}, \extens
{h}{s}{T}) \nonumber\\
&=& \idotsint\Biggl( \prod_{u=t+1}^T \M(x_{u-1},\rmd x_u) g_u(x_u) \Biggr)
h(x_s)\\
&=& \F_{t,T}(x_t, \extens{h}{s}{T}) .\nonumber
\end{eqnarray}
Now, by construction, for any $t \leq s $,
%
%
\begin{equation} \label{eq:two}
\F_{t-1,T}(x_{t-1},\extens{h}{s}{T}) = \int\M(x_{t-1},\rmd x_t)
g_t(x_t) \F_{t,T}(x_t,\extens{h}{s}{T}) .
\end{equation}
The relations \eqref{eq:expressionG}, \eqref{eq:one} and \eqref
{eq:two} imply that
%
%
\begin{equation}
\frac{G^N_{t,T}(x, \extens{h}{s}{T})}{\F_{t,T}(x,\1)} = \frac{\mu
[\F_{t,T}(\cdot, \extens{h}{s}{T})]}{\mu[\F_{t,T}(\cdot,\1)]} -
\frac{\mu'[\F_{t,T}(\cdot,\extens{h}{s}{T})]}{\mu'[ \F
_{t,T}(\cdot,\1)]} , \label{eq:Gn-t-inf-s}
\end{equation}
where $\mu\eqdef\delta_x$ and $\mu'$ is the nonnegative finite
measure defined by
\[
\mu'(A) \eqdef\iint\filt[hat]{t-1}(\rmd x_{t-1}) \M(x_{t-1},\rmd
x_t) g_{t}(x_t)\1_A(x_t) ,\qquad  A \in\Xsigma.
\]
Now, for any finite measure $\mu$ on $(\Xset,\Xsigma)$, the quantity\vspace*{2pt}
\begin{eqnarray*}
&&\frac{\mu[\F_{t,T}(\cdot,\extens{h}{s}{T})]}{\mu[\F_{t,T}(\cdot
,\1)]} \\
&&\qquad = \frac{\idotsint\mu(\rmd x_t) \prod_{u=t+1}^T \M
(x_{u-1},\rmd x_u) g_u(x_{u}) h(x_s)}{\idotsint\mu(\rmd x_t) \prod
_{u=t+1}^T \M(x_{u-1}, \rmd x_u) g_u(x_u)} \\[2pt]
&&\qquad = \frac{\idotsint\mu(\rmd x_t) \prod_{u=t+1}^s \M(x_{u-1},\rmd
x_u) g_u(x_{u}) h(x_s) \F_{s,T}(x_s,\1)}{\idotsint\mu(\rmd x_t)
\prod_{u=t+1}^s \M(x_{u-1},\rmd x_u) g_u(x_u) \F_{s,T}(x_s,\1)}
\end{eqnarray*}
may be seen as the expectation of $h(X_s)$ conditionally on $\chunk
{Y}{t}{T}$, where $X_t$ is distributed according to $A \mapsto\mu(A)/
\mu(\Xset)$. Under the strong mixing condition (Assumption \ref
{assum:strong-mixing-condition}), it is shown in \cite
{delmoralmiclo2000} (see also \cite{delmoral2004}) that, for any\vadjust{\goodbreak} $t
\leq s \leq T$, any finite measure $\mu$ and $\mu'$ on $(\Xset
,\Xsigma)$, any function $h \in\mcb{\Xset}$, that\vspace*{2pt}
\begin{eqnarray*}
&&\biggl| \frac{\idotsint\mu(\rmd x_t) \prod_{u=t+1}^s \M(x_{u-1},\rmd
x_u) g_u (x_u) h(x_s) \F_{s,T}(x_s,\1)}{\idotsint\mu(\rmd x_t) \prod
_{u=t+1}^s \M(x_{u-1},\rmd x_u) g_u (x_u) \F_{s,T}(x_s,\mathbf{1})}
\\[2pt]
&&\quad  {}- \frac{\idotsint\mu'(\rmd x_t) \prod_{u=t+1}^s \M(x_{u-1},\rmd
x_u) g_u (x_u) h(x_s) \F_{s,T}(x_s,\1)}{\idotsint\mu'(\rmd x_t)
\prod_{u=t+1}^s \M(x_{u-1},\rmd x_u) g_u (x_u) \F_{s,T}(x_s,\1)}
\biggr|\\
&&\qquad \leq\rho^{s-t} \oscnorm{h} ,
\end{eqnarray*}
where $\rho$ is defined in \eqref{eq:definition-rho}.
This shows \eqref{eq:time-unif-G} when $t$ is smaller than $s$.

Consider now the case $s < t \leq T$. By definition,\vspace*{2pt}
\begin{eqnarray} \label{eq:cas2premiere}
\F^N_{t,T}(x_t,\extens{h}{s}{T})&=&\idotsint L_{t,T}(\chunk
{x}{0}{t},\extens{h}{s}{T}) \prod_{u=s+1}^{t} \BK{\filt
[hat]{u-1}}(x_u, \rmd x_{u-1}) \nonumber\\[-8pt]\\[-8pt]
&=& \idotsint\F_{t,T}(x_t,\1) \prod_{u=s+1}^{t} \BK{\filt
[hat]{u-1}}(x_u, \rmd x_{u-1}) h(x_s) ,\nonumber
\end{eqnarray}
where the last expression is obtained from the following equality,
valid for $s < t$:\vspace*{2pt}
\begin{eqnarray*}
L_{t,T}(\chunk{x}{0}{t},\extens{h}{s}{T})&=&h(x_s) \idotsint\prod
_{u=t+1}^T \M(x_{u-1},\rmd x_u) g_u(x_u) \\
&=& h(x_s) \F_{t,T}(x_t,\1).
\end{eqnarray*}
Moreover, combining \eqref{eq:cequonveut} and
\eqref{eq:cas2premiere},\vspace*{2pt}
\begin{eqnarray*}
&&\filt[hat]{t-1}[\F^N_{t-1,T}(\cdot,\extens{h}{s}{T})]\\
&&\qquad  = \idotsint
\filt[hat]{t-1}(\rmd u_{t-1}) M(u_{t-1},\rmd x_t) g_t(x_t) \F
^N_{t,T}(x_t,\extens{h}{s}{T}) \\
&&\qquad = \idotsint\filt[hat]{t-1}(\rmd u_{t-1}) M(u_{t-1},\rmd x_t) g_t(x_t)
\F_{t,T}(x_t,\1) \\
&&\hphantom{\qquad = \idotsint}
{}\times\prod_{u=s+1}^t \BK{\filt[hat]{u-1}}(x_u, \rmd
x_{u-1}) h(x_s) .
\end{eqnarray*}
By plugging this expression and \eqref{eq:cas2premiere} into \eqref
{eq:expressionG}, we obtain
\[
\frac{G^N_{t,T}(x, \extens{h}{s}{T})}{\F_{t,T}(x,\1)}=\idotsint\biggl\{
\frac{\mu(\rmd x_{t})}{\mu(\Xset)} - \frac{\mu'(\rmd x_{t})}{\mu
'(\Xset)} \biggr\} \prod_{u=s+1}^t \BK{\filt[hat]{u-1}}(x_u, \rmd
x_{u-1}) h(x_s),
\]
with $\mu(\rmd x_t) = \delta_x(\rmd x_t) \F_{t,T}(x_t,\1)$ and $\mu
'$ being the nonnegative measure defined by
\[
\mu'(A) = \int\filt[hat]{t-1}[\ensuremath{m}(\cdot,x_t)] g_t(x_t)
\F
_{t,T}(x_{t},\1) \1_A(x_{t}) \,\rmd x_t .
\]
Under the uniform ergodicity condition (Assumption \ref
{assum:strong-mixing-condition}) it holds, for any probability measure
$\eta$ on $(\Xset,\Xsigma)$, and any $A \in\Xsigma$,
\[
\BK{\eta}(x,A) = \frac{\int_A \eta( \rmd x') \ensuremath{m}(x',x)
}{\int\eta
( \rmd x') \ensuremath{m}(x',x) } \geq\frac{\sigma_-}{\sigma_+}
\eta(A) ;
\]
thus, the transition kernel $\BK{\eta}$ is uniformly Doeblin with
minorizing constant $\sigma_-/\sigma_+$ and the proof of \eqref
{eq:time-unif-G} for $ s < t \leq T$ follows. The last statement of the
Lemma follows from \eqref{eq:time-unif-G} and the almost-sure convergence
\[
\lim_{N \to\infty}G^N_{t,T}(x,h) = G_{t,T}(x,h) ,\qquad  \PP\mbox{-\as},
\]
for all $x \in\Xset$, which was established in Lemma \ref{lem:limLG}.
\end{pf}

\subsection{A time uniform exponential deviation inequality}
Under the strong mixing Assumption \ref
{assum:strong-mixing-condition}, a time uniform deviation inequality
for the \emph{marginal smoothing} approximation can be derived using
the exponentially decreasing bound on the quantity $G^N_{t,T}$ obtained
in Lemma \ref{lem:G-uniform}.
\begin{thm}
\label{theo:Hoeffding-uniform}
Assume Assumptions \ref{assum:bound-likelihood}--\ref
{assum:strong-mixing-condition} hold with $T= \infty$.
Then, there exist constants $0\leq B,\ C< \infty$ such that for all
integers $N$, $s$, $T$, with $s \leq T$, and for all $\epsilon> 0$,
\begin{eqnarray}
\label{eq:TU-Hoeffding-smoothing-normalized}
\PP[ | \post[hat]{s}{T}(h) - \post{s}{T}(h) | \geq\epsilon] &\leq&
B\rme^{-C N \epsilon^2 / \oscnorm[2]{h}} , \\
\label{eq:TU-Hoeffding-smoothing-normalized-2}
\PP[ | \post[tilde]{s}{T}(h) - \post{s}{T}(h) | \geq\epsilon]
&\leq& B\rme^{-C N \epsilon^2 / \oscnorm[2]{h}} ,
\end{eqnarray}
where $\post[hat]{s}{T}(h)$ and $\post[tilde]{s}{T}(h)$ are defined
in \eqref{eq:smoothing:backw_decomposition_sample} and \eqref
{eq:FFBSi:estimator}.
\end{thm}

Letting $s = T$ in Theorem \ref{theo:Hoeffding-uniform} provides, as a
special case, the (already known) time uniform deviation inequality for
the {\em filter} approximation; however, the novelty of the bounds
obtained here is that these confirm the stability of the FFBSm and
FFBSi marginal smoothing approximations also when $s$ is fixed and $T$
tends to infinity (see \cite{delmoral2004} for further discussion).

\begin{pf*}{Proof of Theorem \ref{theo:Hoeffding-uniform}}
Combining \eqref{eq:Lone} with the definition \eqref{eq:defLt} and
Assumption \ref{assum:strong-mixing-condition} yields, for all $x \in\Xset$,
%
%
\begin{equation} \label{eq:majorationL}
\frac{\sigma_-}{\sigma_+} \leq\frac{\F_{t,T}(x,\1)}{\esssup{\F
_{t,T}(\cdot, \1)}} \leq1 .
\end{equation}
Let $h \in\mcb{\Xset^{T+1}}$ and assume without loss of generality
that $\post{0:T}{T}(h)=0$. Then, \eqref{eq:smooth:recursion} implies
that $\filt{0}[L_{0,T}(\cdot,h)] = 0$ and the first term of the
decomposition \eqref{eq:decomp_Smooth} thus becomes
%
%
\begin{equation}\label{eq:initialRecurInegExpo-1}
\frac{\filt[hat]{0}[L_{0,T}(\cdot,h)]}{\filt[hat]{0}[L_{0,T}(\cdot
,\1)]} = \frac{N^{-1}\sum_{i=0}^N \frac{\rmd\Xinit}{\rmd\XinitIS
{0}}(\epart{0}{i}) g_0(\epart{0}{i}) L_{0,T}(\epart{0}{i},h)}{N^{-1}
\sum_{\ell= 0}^N \frac{\rmd\Xinit}{\rmd\XinitIS{0}}(\epart
{0}{\ell}) g_0(\epart{0}{\ell}) L_{0,T}(\epart{0}{\ell},\1)} ,
\end{equation}
where $(\epart{0}{\ell})_{\ell= 1}^N$ are \iid\ random variables
with distribution $\XinitIS{0}$. Noting that $L_{0,T} = \F_{0,T}$ we
obtain an exponential deviation inequality for \eqref
{eq:initialRecurInegExpo-1} by applying Lemma \ref
{lem:inegessentielle} with
\[
\cases{
\displaystyle a_N = N^{-1}\sum_{i=0}^N \frac{\rmd\Xinit}{\rmd\XinitIS
{0}}(\epart{0}{i}) g_0(\epart{0}{i}) \F_{0,T}(\epart
{0}{i},h)/\supnorm{\F_{0,T}(\cdot,h)} , \cr
\displaystyle b_N = N^{-1}\sum_{i=0}^N \frac{\rmd\Xinit}{\rmd\XinitIS
{0}}(\epart{0}{i}) g_0(\epart{0}{i}) \F_{0,T}(\epart{0}{i},\1
)/\supnorm{\F_{0,T}(\cdot,h)} , \cr
b = \Xinit[g_0(\cdot) \F_{0,T}(\cdot,\1)]/\supnorm{\F_{0,T}(\cdot
,h)} , \cr
\beta= \Xinit(g_0) \sigma_- / \sigma_+ .
}
\]
Here, condition \hyperlink{item:inegEssentielle-borne}{(I)} is trivially
satisfied and conditions \hyperlink{item:inegEssentielle-expo1}{(II)} and \hyperlink{item:inegEssentielle-expo2}{(III)}
follow from the Hoeffding inequality for \iid\ variables.

According to \eqref{eq:decomp_Smooth} and \eqref{eq:definition-A}, it
is now required, for any $1\leq t\leq T$, to derive an exponential
inequality for
\[
A^N_{t,T} \eqdef\frac{N^{-1}\sum_{\ell=1}^N \ewght{t}{\ell}
G^N_{t,T}(\epart{t}{\ell},\extens{h}{s}{T})}{N^{-1} \sum_{\ell
=1}^N \ewght{t}{\ell} \F_{t,T}(\epart{t}{\ell},\1)} .
\]
Note first that, using \eqref{eq:majorationL}, we have
\[
|A^N_{t,T}| \leq\biggl( \frac{\sigma_+}{\sigma_-} \biggr) \frac{N^{-1} \sum
_{\ell=1}^N \ewght{t}{\ell} G^N_{t,T}(\epart{t}{\ell},\extens
{h}{s}{T}) / \esssup{\F_{t,T}(\cdot,\1)}}{N^{-1} \sum_{\ell=1}^N
\ewght{t}{\ell}} .
\]
We use again Lemma \ref{lem:inegessentielle} with
\[
\cases{
\displaystyle a_N = N^{-1} \sum_{\ell=1}^N \ewght{t}{\ell} G^N_{t,T}(\epart
{t}{\ell},\extens{h}{s}{T})/\esssup{\F_{t,T}(\cdot,\1)} , \cr
\displaystyle b_N = N^{-1} \sum_{\ell=1}^N \ewght{t}{\ell} , \cr
b = \mathbb{E}[ \ewght{t}{1} | \mcf{t-1}{N}] = \filt[hat]{t-1} [ \M
(\cdot, g_t) ] / \filt[hat]{t-1}(\adjfunc{t}{t}{}) , \cr
\beta= c_- / \esssup{\adjfunc{t}{t}{}} .
}
\]
Assumption \ref{assum:strong-mixing-condition} shows that $b \geq
\beta$ and Lemma \ref{lem:G-uniform} shows that $|a_N/b_N| \leq M
\eqdef\rho^{|t-s|} \oscnorm{h}$, where $\rho$ is defined in \eqref
{eq:definition-rho}. Therefore, condition~\hyperlink{item:inegEssentielle-borne}{(I)}
of Lemma~\ref{lem:inegessentielle} is
satisfied and the Hoeffding inequality gives
\begin{eqnarray*}
\PP[ | b_N - b | \geq\epsilon] &\leq&\PE\Biggl[ \mathbb{P}\Biggl[ \Biggl|N^{-1} \sum
_{\ell= 1}^N ( \ewght{t}{\ell} - \CPE{\ewght{t}{1}}{\mcf{t-1}{N}}
) \Biggr| \geq\epsilon\Big|\mcf{t-1}{N}\Biggr] \Biggr]\\
&\leq&2 \exp( -2 N \epsilon^2 / \esssup[2]{\ewghtfunc{t}} ),
\end{eqnarray*}
establishing condition \hyperlink{item:inegEssentielle-expo1}{(II)} in Lemma
\ref{lem:inegessentielle}. Finally, Lemma \ref{lem:G-uniform} and the
Hoeffding inequality imply that
\begin{eqnarray*}
\PP[ |a_N | \geq\epsilon] &\leq&\PE\Biggl[ \mathbb{P}\Biggl[ \Biggl|N^{-1} \sum_{\ell
=1}^N \ewght{t}{\ell} G^N_{t,T}(\epart{t}{\ell},\extens
{h}{s}{T})/\esssup{\F_{t,T}(\cdot,\1) } \Biggr| \geq\epsilon\Big|\mcf
{t-1}{N}\Biggr] \Biggr]\\
&\leq&2 \exp\biggl( - 2 \frac{N \epsilon^2}{\esssup[2]{\ewghtfunc{t}}
\rho^{2|t-s|} \oscnorm[2]{h}} \biggr) = 2 \exp\biggl( - 2 \frac{N \epsilon
^2}{\esssup[2]{\ewghtfunc{t}} M^2} \biggr) .
\end{eqnarray*}
Lemma \ref{lem:inegessentielle} therefore yields
\[
\PP\biggl( \biggl| \frac{a_N}{b_N} \biggr| \geq\epsilon\biggr) \leq2 \exp\biggl( - \frac{N
\epsilon^2 c_-^2}{2 \oscnorm[2]{h} \rho^{2|t-s|}\esssup{\ewghtfunc
{t}}^2\esssup{\adjfunc{t}{t}{}}^2} \biggr) ,
\]
so that
\[
\PP( | A^N_{t,T} | \geq\epsilon) \leq2 \exp\biggl( - \frac{N \epsilon
^2 c_-^2 \sigma_-^2}{2 \oscnorm[2]{h} \rho^{2|t-s|}\esssup
{\ewghtfunc{t}}^2\esssup{\adjfunc{t}{t}{}}^2 \sigma_+^2} \biggr) .
\]
A time uniform exponential deviation inequality for $\sum_{t=1}^T
A_{t,T}$ then follows from Lemma \ref{lem:pasfor} and the proof is complete.
\end{pf*}

\subsection{A time uniform bound on the variance of the marginal
smoothing distribution}
\label{sec:TimeUniformCLTFFBS}
Analogous to the result obtained in the previous section, a~time
uniform bound on the asymptotic variance in the CLT for the {\em
marginal smoothing} approximations can, again under the strong mixing
Assumption~\ref{assum:strong-mixing-condition}, be easily obtained
from the exponentially decreasing bound on $G_{t,T}$ stated and proved
in Lemma \ref{lem:G-uniform} for the quantity.

\begin{thm} \label{theo:CLT-uniform}
Assume Assumptions \ref{assum:bound-likelihood}--\ref
{assum:strong-mixing-condition} hold with $T = \infty$. Then, for all
${s\leq T}$,
\[
\asymVar{0:T}{T}{\extens{h}{s}{T}} \leq\biggl(\frac{\sigma_+}{\sigma
_-} \Bigl( 1 \vee\sup_{t \geq1} \supnorm{\adjfunc{t}{t}{}} \Bigr) \sup_{t
\geq0} \supnorm{\ewghtfunc{t}}\oscnorm{h} \biggr)^2 \frac{1+\rho
^2}{1-\rho^2} ,
\]
where $\asymVar{0:T}{T}{}$ is defined in \eqref{eq:expression-covariance}.
\end{thm}

In accordance with the results of the previous section, letting $s = T$
in the previous theorem provides a time uniform bound on the asymptotic
variance for the {\em filter} approximation; nevertheless, as mentioned
previously, the situation of interest for us is when $s$ is fixed and
$T$ goes to infinity.
\begin{pf*}{Proof of Theorem \ref{theo:CLT-uniform}}
Combining \eqref{eq:majorationL} and \eqref{eq:time-unif-G-lim} with
$\XinitIS{0}(\ewghtfunc{0})=1$ yields
\[
\frac{\XinitIS{0}(\ewghtfunc{0}^2(\cdot) G_{0,T}^2(\cdot,\extens
{h}{s}{T}))}{\XinitIS{0}^2 [\ewghtfunc{0}(\cdot) \F_{0,T}(\cdot,\1
) ]} \leq\biggl(\frac{\sigma_+}{\sigma_-}\supnorm{\ewghtfunc
{0}}\oscnorm{h}\rho^s \biggr)^2 .
\]
Moreover, by inserting, for any $0 < t \leq T$, the bound obtained in
\eqref{eq:time-unif-G-lim} into the expression \eqref
{eq:definition-upsilon} of $\upsilon_{t,T}$ we obtain
\[
\frac{\filt{t-1}(\upsilon_{t,T}(\cdot,\extens{h}{s}{T}))) \filt
{t-1}(\adjfunc{t}{t}{})}{\filt{t-1}^2[\F_{t-1,T}(\cdot,\1)]} \leq
\biggl(\frac{\sigma_+}{\sigma_-} \supnorm{\adjfunc{t}{t}{}}\supnorm
{\ewghtfunc{t}}\oscnorm{h} \rho^{|t-s|} \biggr)^2 .
\]
Finally, plugging the two bounds above into \eqref
{eq:expression-covariance} gives
\[
\asymVar{0:T}{T}{\extens{h}{s}{T}}\leq\biggl(\frac{\sigma_+}{\sigma_-}
\Bigl( 1 \vee\sup_{t \geq1} \supnorm{\adjfunc{t}{t}{}} \Bigr) \sup_{t \geq
0} \supnorm{\ewghtfunc{t}}\oscnorm{h} \biggr)^2 \Biggl(\sum_{t=0}^\infty\rho
^{2|t-s|} \Biggr) ,
\]
which completes the proof.
\end{pf*}

\section{\texorpdfstring{Proofs of Propositions \protect\ref{prop:complexityBoostrap} and \protect\ref{prop:complexity}}
{Proofs of Propositions 1 and 2}}
\label{sec:complexity:proofs}
Having at hand the theory established in the previous sections, we are
now ready to present the proofs of Propositions \ref
{prop:complexityBoostrap} and~\ref{prop:complexity}.

\begin{pf*}{Proof of Proposition \ref{prop:complexityBoostrap}}
The average number of simulations\vspace*{-2pt} required to sample $J_s^\ell$
conditionally on $\mathcal{G}_{s+1}^N$ is $\sigma_+ \sumwght{s} /
\sum_{u = 1}^N \ewght{s}{u} \ensuremath{m}(\epart{s}{u}, \epart
{s+1}{J_{s+1}^\ell})$. Hence, the number of simulations $Z_s^N$
required to sample $\{ J_s^\ell\}_{\ell= 1}^N$ has conditional expectation
\[
\mathbb{E} [ Z_s^N | \mathcal{G}_{s+1}^N ] = \sum_{\ell= 1}^N \frac
{\sigma_+ \sumwght{s}}{\sum_{i = 1}^N \ewght{s}{i} \ensuremath
{m}(\epart
{s}{i},\epart{s+1}{J_{s+1}^\ell})} .
\]
We denote $\ewght{s|T}{i} \eqdef\mathbb{P}[J^1_s = i|\mcf{T}{N}]$
and $\ewght{s:s+1|T}{\ell i} \eqdef\mathbb{P}[J^1_s = \ell,
J^1_{s+1} = i|\mcf{T}{N}]$ and write
\begin{eqnarray*}
\mathbb{E} [ Z_s^N | \mcf{T}{N} ] &=& \sum_{i=1}^N \ewght
{s+1|T}{i}\frac{\sigma_+ \sumwght{s}}{\sum_{j=1}^N \ewght{s}{j}
\ensuremath{m}
(\epart{s}{j}, \epart{s+1}{i})} \\
&=& \sigma_+ \sumwght{s} \sum_{i=1}^N \sum_{\ell= 1}^N \frac
{\ewght{s+1|T}{i} \ewght{s}{\ell} \ensuremath{m}(\epart{s}{\ell},
\epart
{s+1}{i})}{\sum_{j=1}^N \ewght{s}{j} \ensuremath{m}(\epart{s}{j},
\epart
{s+1}{i})} \times\frac{1}{\ewght{s}{\ell} \ensuremath{m}(\epart
{s}{\ell
},\epart{s+1}{i})} \\
&=& \sigma_+ \sumwght{s} \sum_{i=1}^N \sum_{\ell= 1}^N \ewght
{s:s+1|T}{\ell i} \frac{1}{\ewght{s}{\ell} m(\epart{s}{\ell},
\epart{s+1}{i})} .
\end{eqnarray*}
For the bootstrap particle filter, $\ewght{s}{\ell} \equiv g_s(\epart
{s}{\ell})$; Theorem \ref{thm:Hoeffding-FFBS} then implies that
$\sumwght{s}/ N \plim_{N \rightarrow\infty} \post{s}{s-1}(g_s)$ and
\begin{eqnarray*}
&&\sum_{i=1}^N \sum_{\ell= 1}^N \ewght{s:s+1|T}{\ell i} \frac
{1}{\ewght{s}{\ell} \ensuremath{m}(\epart{s}{\ell}, \epart
{s+1}{i})} \\
&&\qquad \plim
_{N \rightarrow\infty} \iint\post{s:s+1}{T}(\rmd x_{s:s + 1}) \frac
{1}{g_s(x_s) m(x_s, x_{s+1})} .
\end{eqnarray*}
Besides,
\begin{eqnarray*}
&&\iint\post{s:s+1}{T}(\rmd x_{s:s + 1}) \frac{1}{g_s(x_s)
\ensuremath{m}(x_s, x_{s+1})} \\
&&\qquad  = \idotsint\post{s}{s-1}(\rmd x_s) \frac{g_s(x_s) \M(x_s,
\rmd x_{s+1})}{g_s(x_s) \ensuremath{m}(x_s, x_{s+1})}
g_{s+1}(x_{s+1})\\
&&\qquad \quad \hphantom{\idotsint}
{}\times\prod_{u=s+2}^T \M(x_{u-1}, \rmd x_u) g_u(x_u)\\
&&\qquad \quad {}\Big/\idotsint\post
{s}{s-1}(\rmd x_s) g_{s}(x_{s}) \prod_{u=s+1}^T \M(x_{u-1}, \rmd x_u)
g_u(x_u) \\
&&\qquad  = \frac{\idotsint\rmd x_{s+1} \prod_{u=s+2}^T \int\M(x_{u-1},
\rmd x_u) g_u(x_u)}{\idotsint\post{s}{s-1}(\rmd x_s) g_s(x_s) \prod
_{u=s+1}^T \M(x_{u-1}, \rmd x_u) g_u(x_u)} .
\end{eqnarray*}
Similarly, in the fully adapted case we have $\ewght{s}{i} \equiv1$
for all $i \in\{1, \dots, N\}$; thus, $\sumwght{s} = N$ and
\begin{eqnarray*}
&&\sum_{i=1}^N \sum_{\ell= 1}^N \ewght{s:s+1|T}{\ell i} \frac
{1}{\ewght{s}{\ell} \ensuremath{m}(\epart{s}{\ell}, \epart
{s+1}{i})} \\
&&\qquad \plim
_{N \rightarrow\infty} \iint\post{s:s+1}{T}(\rmd x_{s:s + 1}) \frac
{1}{m(x_s, x_{s+1})}
\\
&&\qquad = \frac{\idotsint g_{s+1}(x_{s+1})\, \rmd x_{s+1}\prod_{u=s+2}^T
\int\M(x_{u-1}, \rmd x_u) g_u(x_u)}{\idotsint\filt{s}(\rmd x_s)
\prod_{u=s+1}^T \M(x_{u-1}, \rmd x_u) g_u(x_u)} .
\end{eqnarray*}
In both cases, the numerator can be bounded from above by
\[
\frac{\sigma_+^{T-s-1} \prod_{u=s+1}^T \int g_u(x_u)\, \rmd
x_u}{\idotsint\post{s}{s-1}(\rmd x_s) g_s(x_s) \prod_{u=s+1}^T \M
(x_{u-1}, \rmd x_u) g_u(x_u)}
\]
if $\int g_u(x_u) \,\rmd x_u < \infty$ for all $u \geq0$.
\end{pf*}

\begin{pf*}{Proof of Proposition \ref{prop:complexity}}
Fix a time step $s$ of the algorithm and denote by $C_s$ the number of
elementary operations required for this step. For $k \in\{1, \dots,
n\}$, let $T_s^k$ be the number
of times that $k$ appears in list $L$ at time $s$\vspace*{-2pt} in the `while' loop.
Let also $N_s^u \eqdef\sum_{k=1}^N \mathbh{1}_{\{T_s^k \geq u\}}$
be the size of~$L$ (i.e., the value of $n$ at line 6)\vspace*{-4pt} after $u$
iterations of the `while' loop, with $N_s^0 \eqdef N$. Then, using
Proposition \ref{prop:multisample} there exists a constant C such that
\[
C_s \leq
C \sum_{u=0}^\infty N_s^u \biggl( 1 + \log\biggl(1 + \frac{N}{N_s^u} \biggr) \biggr) .
\]
As $n\to n(1+\log(1+N/n))$ is a concave, increasing function, it holds
by Jensen's inequality that
\[
\PE[C_s] \leq C \sum_{u=0}^\infty\PE[N_s^u] \biggl( 1+\log\biggl(1 + \frac
{N}{\PE[N_s^u]} \biggr) \biggr) .
\]
Besides,
\[
\PE[N_s^u] = \sum_{k = 1}^N \PP(T_s^k \geq u) \leq N \biggl(1 - \frac
{\sigma_-}{\sigma_+} \biggr)^u
\]
as $\sigma_-/\sigma_+$ is a lower bound on the acceptation probability.
Thus,
\[
\PE[C_s] \leq C N \sum_{u=0}^\infty\biggl(1-\frac{\sigma_-}{\sigma_+}
\biggr)^{u} \biggl( 1+\log\biggl(1+\frac{1}{ (1-\sigma_-/\sigma_+ )^u} \biggr) \biggr)\leq
\frac{KN\sigma_+}{\sigma_-} .
\]
\upqed
\end{pf*}

\begin{appendix}
\section{\texorpdfstring{Proof of Lemma \lowercase{\protect\ref{lem:inegessentielle}}}
{Appendix A: Proof of Lemma 4}}\label{sec:proof:lem:inegEssentielle}
Write
\[
\biggl|\frac{a_N}{b_N} \biggr| \leq b^{-1} \biggl| \frac{a_N}{b_N} \biggr| |b-b_N| + b^{-1} |
a_N | \leq\beta^{-1}M |b-b_N|+ \beta^{-1} | a_N | ,\qquad
\mathbb{P}\mbox{-a.s.}
\]
Thus,
\[
\biggl\{ \biggl|\frac{a_N}{b_N} \biggr| \geq\epsilon\biggr\} \subseteq\biggl\{|b-b_N| \geq\frac
{\epsilon\beta}{2M} \biggr\} \cup\biggl\{ |a_N | \geq\frac{\epsilon\beta
}{2} \biggr\} ,
\]
from which the proof follows.

\section{\texorpdfstring{Technical results}{Appendix B: Technical results}}\label{sec:TechnicalResults}

\begin{lem} \label{lem:pasfor}
Let $\{ Y_{n,i} \}_{i=1}^n$ be a triangular array of random variables
such that there exist constants $B > 0$, $C > 0$, and $\rho$ with $0 <
\rho< 1$ satisfying, for all $n$, $i \in\{1, \dots, n \}$, and
$\epsilon> 0$,
\[
\PP( |Y_{n,i}| \geq\epsilon) \leq B \exp( -C \epsilon^2 \rho
^{-2i} ) .
\]
Then, there exist constants $\bar{B} > 0$ and $\bar{C} > 0$ such
that, for any $n$ and $\epsilon> 0$,
\[
\PP\Biggl(\Biggl | \sum_{i=1}^n Y_{n,i} \Biggr| \geq\epsilon\Biggr) \leq\bar{B} \rme
^{- \bar{C} \epsilon^2} .
\]
\end{lem}

\begin{pf}
Set $S \eqdef\sum_{\ell= 1}^\infty\sqrt{\ell} \rho^\ell$; one
easily concludes that
\[
\PP\Biggl( \Biggl| \sum_{i=1}^n Y_{n,i} \Biggr| \geq\epsilon\Biggr) \leq\sum_{i=1}^n \PP
\bigl( |Y_{n,i}| \geq\epsilon S^{-1} \sqrt{i} \rho^i \bigr) \leq B \sum
_{i=1}^n \exp( -C S^{-1} \epsilon^2 i ) .
\]
Set $\epsilon_0 > 0$. The proof follows by noting that, for any
$\epsilon\geq\epsilon_0$,
\[
\sum_{i=1}^n \exp( -C S^{-1} \epsilon^2 i ) \leq\exp( C S^{-1}
\epsilon_0^2 ) \exp( -C S^{-1} \epsilon^2 ) / \bigl( 1 - \exp( C S^{-1}
\epsilon_0^2 ) \bigr) .\qquad
\]
\upqed
\end{pf}

\subsection{Description of the sampling procedure}\label{subsec:multisampling}
In this section, we describe and analyze an efficient multinomial
sampling procedure, detailed in Algorithm~\ref{alg:multisample}. Given
a probability distribution $(p_1,\dots,p_N)$ on the set $\{1, \dots,
N\}$, it returns a sample of size $n$ of that distribution. Compared to
the procedure described in Section 7.4.1 in \cite
{cappemoulinesryden2005}, its main virtue is to be efficient for both
large and small samples sizes: if $n = 1$, the complexity is $O(\log
(N))$, while if $n = N$, the complexity is $O(N)$.

\begin{prop} \label{prop:multisample}
The number of elementary operations required by Algorithm \ref
{alg:multisample} is $O (n + n \log(1 + N/n) )$.
\end{prop}

\begin{pf}
The order statistics at line 5 and the permutation at line 6 can be
sampled using $O(n)$ operations; see \cite{devroye1986}, Chapter V and XIII.
For each value of $k$ between $1$ and $n$, denote by $G_k$ the number
of times lines 11--13 are executed. Observe that line 18 is executed
the same number of times, and thus the number of elementary operations
required by call to Algorithm~\ref{alg:multisample} is $O(n + \sum
_{k=1}^nG_k)$. But the value of $l$ is increased during iteration $k$
by at least $2^{G_k}-1$, and as the final value of $l$ is at most equal
to $N$, it holds that
\[
\sum_{k=1}^n 2^{G_k} \leq N + n .
\]
By convexity,
\[
\exp\Biggl(\frac{\log(2)}{n} \sum_{k=1}^n G_k \Biggr) \leq\frac{1}{n} \sum
_{k=1}^n 2^{G_k} 
\leq1 + \frac{N}{n} ,
\]
which implies that
\[
\sum_{k=1}^nG_k \leq
n \log\biggl(1 + \frac{N}{n} \biggr) / \log(2) .
\]
\upqed
\end{pf}\vfill\eject

\begin{algorithm}[t]
\begin{algorithmic}[1]
\caption{Multinomial sampling} \label{alg:multisample}
\State$q_1 \gets p_1$
\For{$k$ from 1 to $N$}
\State$q_k \gets q_{k-1}+p_k$
\EndFor
\State sample an order statistics $U_{(1)},\dots,U_{(n)}$ of an i.i.d.
uniform distribution
\State uniformly sample a permutation $\sigma$ on $\{1, \dots, n\}$
\State$l \gets0, r\gets1$
\For{$k$ from $1$ to $n$}
\State$d \gets1$
\While{$U_{(k)} \geq q_r$}
\State$l \gets r$
\State$r \gets\min(r + 2^d, N)$
\State$d \gets d+1$
\EndWhile
\While{$r-l>1$}
\State$m \gets\lfloor(l+r)/2 \rfloor$
\If{$U_{(k)}\geq q_m$}
\State$l \gets m$
\Else
\State$r \gets m$
\EndIf
\EndWhile
\State$I_{\sigma(k)} \gets r$
\EndFor
\end{algorithmic}
\end{algorithm}

\end{appendix}

%

\printaddresses

\end{document}